\documentclass[11pt]{amsart}

\usepackage{verbatim}

\usepackage{anysize}

\marginsize{3cm}{3cm}{2cm}{2cm}

\usepackage{epsf}

\usepackage{euscript}

\usepackage{pb-diagram}

\usepackage{amssymb}

\def\to{\mathchoice
{\longrightarrow}
{\rightarrow}
{\rightarrow}
{\rightarrow}}

\def\mapsto{\mathchoice
{\DOTSB\mapstochar\longrightarrow}
{\DOTSB\mapstochar\rightarrow}
{\DOTSB\mapstochar\rightarrow}
{\DOTSB\mapstochar\rightarrow}}

\def\hookrightarrow{\mathchoice
{\DOTSB\lhook\joinrel\relbar\joinrel\rightarrow}
{\DOTSB\lhook\joinrel\rightarrow}
{\DOTSB\lhook\joinrel\rightarrow}
{\DOTSB\lhook\joinrel\rightarrow}}

\def\subsection#1{\refstepcounter{subsection}
\medskip\noindent
{\textbf{(\thesubsection)\ #1\unskip. }}\ignorespaces}

\def\subsubsection#1{\refstepcounter{subsubsection}
\smallskip\noindent
{\textbf{(\thesubsubsection)}\ \textit{#1\unskip. }}\ignorespaces}

\def\theequation{\thesection.\arabic{equation}}
\def\thesubsection{\thesection.\arabic{subsection}}
\def\thesubsubsection{\thesection.\arabic{subsection}.\arabic{subsubsection}}

\makeatletter

\def\@setcopyright{}
\def\serieslogo@{}

\let\c@equation=\c@subsection

\def\ref#1{\@ifundefined{r@#1}%
         {{\bf ??}%
          \iffirstchoice@
             \@warning
                {Reference `#1' on page \thepage\space undefined}%
          \fi}%
         {(\expandafter\expandafter\expandafter
          \@car \csname r@#1\endcsname \@nil\null)}}

\makeatother

\newtheorem{theorem}[subsection]{Theorem}
\newtheorem{proposition}[subsection]{Proposition}
\newtheorem{lemma}[subsection]{Lemma}
\newtheorem{corollary}[subsection]{Corollary}

\theoremstyle{definition}
\newtheorem{definition}[subsection]{Definition}

\newtheorem{remark}[subsection]{Remark}

\def\Z{\mathbb{Z}}
\def\C{\mathbb{C}}

\def\R{\mathbb{R}}
\def\Q{\mathbb{Q}}

\newcommand{\om}{\omega}
\renewcommand{\o}{\otimes}

\DeclareMathOperator{\Hom}{Hom}
\DeclareMathOperator{\Aut}{Aut}
\DeclareMathOperator{\Ind}{Ind}
\DeclareMathOperator{\Res}{Res}

\DeclareMathOperator{\Iso}{Iso}

\DeclareMathOperator{\gr}{gr}
\DeclareMathOperator{\ch}{ch}
\DeclareMathOperator{\Ch}{Ch}
\DeclareMathOperator{\CCh}{\mathbb{C}\textup{h}}
\DeclareMathOperator{\rk}{rk}
\DeclareMathOperator{\Tr}{Tr}
\DeclareMathOperator{\Exp}{Exp}
\DeclareMathOperator{\Log}{Log}
\DeclareMathOperator{\sgn}{\varepsilon}
\DeclareMathOperator{\VV}{Vert}
\DeclareMathOperator{\Edge}{Edge}
\DeclareMathOperator{\Leg}{Leg}
\DeclareMathOperator{\Flag}{Flag}
\DeclareMathOperator{\Ob}{Ob}
\DeclareMathOperator{\Mor}{Mor}
\DeclareMathOperator{\Id}{Id}
\DeclareMathOperator*{\colim}{colim}

\DeclareMathOperator{\OR}{\operatorname{Or}}
\DeclareMathOperator{\Cyc}{Cyc}
\DeclareMathOperator{\Spec}{Spec}
\DeclareMathOperator{\sd}{sd}

\DeclareMathOperator{\Det}{Det}
\def\local{\mathfrak{D}}

\newcommand{\alg}{{\text{alg}}}

\DeclareMathOperator{\Ass}{\mathcal{A}\mathit{ss}}
\DeclareMathOperator{\Com}{\mathcal{C}\!\mathit{om}}
\DeclareMathOperator{\Lie}{\mathcal{L}\mathit{ie}}

\renewcommand{\a}{{\EuScript A}}
\renewcommand{\b}{{\EuScript B}}
\renewcommand{\c}{{\EuScript C}}
\newcommand{\e}{{\EuScript E}}
\newcommand{\CE}{{\EuScript E}}
\renewcommand{\v}{{\EuScript V}}
\newcommand{\w}{{\EuScript W}}
\newcommand{\x}{{\EuScript X}}
\renewcommand{\SS}{\mathbb{S}}
\newcommand{\bull}{{\bullet}}
\newcommand{\CS}{{\EuScript S}}

\newcommand{\CC}{{\EuScript C}}

\newcommand{\CG}{{\EuScript G}}
\newcommand{\CU}{{\EuScript U}}
\newcommand{\Cat}{{\EuScript C}}
\newcommand{\1}{{1\!\!1}}
\let\eps=\varepsilon
\let\phi=\varphi
\newcommand{\op}{{\textrm{op}}}
\newcommand{\TOP}{{\textrm{top}}}
\newcommand{\p}{\partial}

\newcommand{\CM}{{\EuScript M}}
\newcommand{\Mhat}{\widehat{\EuScript M}}
\newcommand{\Cinf}{C^\infty}
\newcommand{\CL}{{\EuScript L}}
\newcommand{\half}{{\textstyle\frac12}}

\newcommand{\pp}{{\EuScript P}}
\renewcommand{\*}{\cdot}
\renewcommand{\L}{{\mathsf L}}
\newcommand{\RR}{{\mathsf R}}
\newcommand{\TT}{\mathbb{T}}

\def\<{\langle}
\def\>{\rangle}

\def\k{{\mathbf k}}

\newcommand{\Mbar}{\overline{\EuScript M}}
\newcommand{\tGrav}{{\mathsf{Grav}}}
\newcommand{\Grav}{{\mathcal{G}\mathit{rav}}}

\def\.{\wedge}
\def\({(\!(}
\def\){)\!)}

\def\]{{]\!]}}
\def\[{{[\![}}
\def\MM{\mathbb{M}}

\def\FF{{\mathsf F}}
\def\BB{{\mathsf B}}
\newcommand{\Coprod}{{\textstyle\coprod}}
\def\CP{\mathbb{CP} }
\newcommand{\Dual}{\vee}
\newcommand{\tom}{\tilde\om}

\newcommand{\twist}{{\mathfrak{T}}}
\newcommand{\can}{{\mathfrak{K}}}
\DeclareMathOperator{\pv}{pv}

\newcommand{\tildetilde}{\hat}

\newcommand{\nf}{\mathcal{NF}}
\newcommand{\df}{\mathcal{DF}}

\newcommand{\+}{+{}}

\newcommand{\Susp}{{\mathfrak{s}}}
\newcommand{\lam}{{\mathfrak{l}}}
\newcommand{\Pp}{{\mathfrak{p}}}

\begin{document}

\title{Modular operads}

\author{E. Getzler and M.M. Kapranov}

\address{Department of Mathematics, MIT, Cambridge, Massachusetts 02139
USA}

\email{getzler@math.mit.edu}

\address{Department of Mathematics, Northwestern University, Illinois
60208 USA}

\email{kapranov@math.nwu.edu}


\maketitle

\section*{Introduction}

Recently, there has been increased interest in applications of operads
outside homotopy theory, much of it due to the relation between operads and
moduli spaces of algebraic curves.

The formalism of operads is closely related to the combinatorics of trees
\cite{GJ}, \cite{GK}. However, in dealing with moduli spaces of curves, one
encounters general graphs, the case of trees corresponding to curves of
genus $0$.

This suggests considering a ``higher genus'' analogue of the theory of
operads, in which graphs replace trees. We call the resulting objects
modular operads: their systematic study is the purpose of this paper.

The cobar functor $\BB$ \cite{GK} is an involution on the category of
differential graded (dg) operads. We will construct an analogous functor
$\FF$ on the category of dg-modular operads, the Feynman transform. This
functor generalizes Kontsevich's graph complexes \cite{Kontsevich}.

The behaviour of $\FF$ is more mysterious than that of the cobar
construction. For example, for such a simple operad as $\Com$, describing
commutative algebras, $\BB\Com$ is a resolution of the Lie operad. On the
other hand, knowledge of the homology of $\FF\Com$ implies complete
information on the dimensions of the spaces of Vassiliev invariants of
knots (by a theorem of Kontsevich and Bar-Natan \cite{BN}; see
\ref{Vassiliev}).

Our main result about the Feynman transform is the calculation of its Euler
characteristic; to do this, we use the theory of symmetric functions. As a
model for this calculation, take the formula for the enumeration of graphs
known in mathematical physics as Wick's theorem \cite{BIZ}. Consider the
asymptotic expansion of the integral
\begin{equation} \label{w}
W(\xi,\hbar) = \log \int \exp \frac{1}{\hbar} \Bigl( x\xi - \frac{x^2}{2} +
\sum_{2(g-1)+n>0} \frac{a_{g,n}\hbar^gx^n}{n!}  \Bigr) \,
\frac{dx}{\sqrt{2\pi\hbar}} ,
\end{equation}
considered as a power series in $\xi$ and $\hbar$. (The asymptotic
expansion is independent of the domain of integration, provided it contains
$0$.) Let $\Gamma\(g,n\)$ be the set of isomorphism classes of connected
graphs $G$, with a map $v\mapsto g(v)$ from the vertices $\VV(G)$ of $G$ to
$\{0,1,2,\dots\}$ and having exactly $n$ legs numbered from $1$ to $n$,
such that
$$
g = \sum_{v\in\VV(G)} g(v) + b_1(G) ,
$$
where $b_1(G)$ is the first Betti number of the graph. If $v$ is a vertex
of $G$, denote by $n(v)$ its valence, and let $|\Aut(G)|$ be the
cardinality of the automorphism group of $G$. Wick's theorem states that
\begin{equation} \label{wick}
W \sim \frac{1}{\hbar} \Bigl( \frac{\xi^2}{2} + \sum_{2(g-1)+n>0}
\frac{\hbar^g \xi^n}{n!}
\sum_{G\in\Gamma\(g,n\)} \frac{1}{|\Aut(G)|} \prod_{v\in\VV(G)}
a_{g(v),n(v)} \Bigr) .
\end{equation}

The calculation of the Euler characteristic of $\FF$ is a natural
generalization of this, in which the coefficients $a_{g,n}$ are replaced by
representations $\v\(g,n\)$ of the symmetric groups $\SS_n$, sums and
products are replaced by the operations $\oplus$ and $\o$, and the weight
$|\Aut(G)|^{-1}$ is replaced by taking the coinvariants with respect to a
natural action of $\Aut(G)$. Up to isomorphism, a sequence
$\v=\{\v\(g,n\)\mid 2(g-1)+n>0\}$ of $\SS_n$-modules is determined by its
Frobenius characteristic, which is a symmetric function $f(x_1,x_2,\dots)$
in infinitely many variables. In Sections~7 and 8, we define analogues of
the Legendre and Fourier transforms for symmetric functions, which give
formulas for the characteristics of $\BB\a$ and $\FF\a$, where $\a$ is a
cyclic, respectively modular, operad.

There is an Euler characteristic associated to orbifolds, called the
virtual Euler characteristic. This is the invariant obtained by giving each
cell $\sigma$ a coefficient $(-1)^{\dim(\sigma)}/|\Aut(\sigma)|$, very much
like in \ref{wick}. Harer and Zagier \cite{HZ} calculated the virtual Euler
characteristics of the moduli spaces $\CM_{\gamma,\nu}$ and
$\CM_{\gamma,\nu}/\SS_\nu$ (see also Penner \cite{Penner:euler}), and more
recently, Kontsevich has given a very simple proof of their result using
Wick's theorem \cite{Kontsevich:KdV}.

Let $|\CM_{\gamma,\nu}|$ be the coarse moduli space of smooth algebraic
curves of genus $\gamma$ with $\nu$ labelled marked points, and
$|\CM_{\gamma,\nu}|/\SS_\nu$ that of smooth algebraic curves of genus
$\gamma$ with $\nu$ unlabelled marked points. Let $e(|\CM_{\gamma,\nu}|)$
and $e(|\CM_{\gamma,\nu}|/\SS_\nu)$ be their Euler characteristics;
clearly, these coincide for $\nu=1$. Harer and Zagier were able to
calculate $e(|\CM_{\gamma,1}|)$; however, for higher values of $\nu$,
little is known about Euler characteristics $e(|\CM_{\gamma,\nu}|)$ and
$e(|\CM_{\gamma,\nu}|/\SS_\nu)$. By applying our formulas to the modular
operad $\Ass$ corresponding to associative algebras, we obtain in Section~9
a closed formula for the sum
$$
\sum_{2(1-\gamma)-\nu=\chi} e(|\CM_{\gamma,\nu}|/\SS_\nu) ,
$$
where $\chi$ is a fixed integer, representing the Euler characteristic of
the punctured Riemann surfaces contributing to the sum.

The use of symmetric functions in enumeration of graphs goes back to
P\'olya \cite{Polya}. Our approach is slightly different: while he
associates symmetric functions to permutations of vertices of the graph, we
associate them to permutations of flags of the graph (pairs consisting of a
vertex and an incident edge). The idea of attaching arbitrary
representations of symmetric groups to vertices of a tree appears (under
the name ``lumps'') in Hanlon-Robinson \cite{HR}; they obtain formulas
resembling our formula for the characteristic of $\BB\a$ (in P\'olya's
setting). The introduction of the Legendre transform in this problem leads
to a new perspective on this class of problems by bringing out a hidden
involutive symmetry, which is very natural from the point of view of
operads.

Our analogue of Wick's theorem may be viewed as a synthesis of the methods
of graphical enumeration of quantum field theory with P\'olya's ideas. Our
formula for the character of $\FF\a$ has another link to quantum field
theory, since the space of symmetric functions is the Hilbert space for the
basic representation of $\operatorname{GL}_{\text{res}}(\infty)$ (Kac-Raina
\cite{KR}); in this direction, we present a formal representation of the
characteristic of the free modular operad $\MM\v$ as a functional integral
\ref{Wick}.

We now describe the contents of this paper. In Section 1, we recall the
definition of cyclic operads from \cite{cyclic}; roughly speaking, these
are operads in which the inputs and output may be permuted. In Section 2,
we define modular operads as algebras over a certain triple, constructed by
summing over graphs. Modular operads are actually a special sort of cyclic
operad, in which there is an additional operation (which we call
contraction) which reduces the number of inputs by two. In Section 3, we
explain the axioms which must be imposed on such a contraction in order
that it determines a modular operad structure; this may be viewed as a
coherence theorem for modular operads.

In Section 4, we introduce a generalization of modular operads, in which
certain signs are introduced into the structure maps, which we call
cocycles. A cocycle is a certain functor from graphs to the Picard category
of invertible graded vector spaces. The most important cocycle for us will
be the determinant of the first cohomology of the graph, which is obviously
trivial when restricted to trees; this explains why this twist is not
needed in the theory of operads. In Section 5, we construct the Feynman
transform $\FF$, which maps from the category of dg-modular operads to the
category of dg-modular operads for this cocycle.

Another example of a Feynman transform is given in Section~6: roughly
speaking, the complexes of currents on the moduli spaces of stable curves
$\Mbar_{g,n}$ form a modular operad, and the Feynman transform of this
modular operad may be identified with the differential forms on the open
strata $\CM_{g,n}$; these form a twisted modular operad, in the sense of
Section 4, by a construction involving residues taken around divisors at
infinity.

In Sections~7 and 8, we define analogues of the Legendre and Fourier
transforms for symmetric functions. In this way, we obtain formulas for the
characteristics of $\BB\a$ and $\FF\a$, where $\a$ is a cyclic,
respectively modular, operad.

One of the pioneers of the use of Wick's theorem as a tool in topology was
Claude Itzykson, and he was an influence on us and many of our colleagues,
in innumerable ways. We humbly offer this article in his memory.

\section*{Acknowledgements}

First of all, we must thank P. Hanlon for sharing his ideas on the role of
symmetric functions in the theory of operads. To a great extent, our
interest in the subject was stimulated by talks of and conversations with
M. Kontsevich. We would like to thank D. Zagier for a very useful
discussion on our results. We must also thank D. Bar-Natan and R. Bott for
informing us of an error (in the definition of the Feynman transform) in a
previous version of this paper.

The research of both authors was partially supported by the NSF and the
A.P. Sloan Foundation. This paper was written while the authors were guests
of the Max Planck Institute for Mathematics in Bonn, whose hospitality and
financial support is greatly appreciated.

\section{Cyclic operads}

In this section, we recall the definition of a cyclic operad --- this will
be useful later, since one way of looking at modular operads is as a
special kind of cyclic operad.

Our presentation of the theory of cyclic operads is a little different from
our previous account \cite{cyclic}; we need a non-unital version of the
theory, due to Markl \cite{Markl}. One advantage of this formulation is
that the basic operations in an operad are bilinear. In any case, if one
simply took the original definition of an operad (May \cite{May}), and
omitted the axioms involving the unit, one would not obtain the same
notion.

\subsection{$\SS$-modules}
Throughout this paper, we work over a fixed field $\k$ of characteristic
$0$.

A chain complex (dg-vector space) is a graded vector space $V_\bull$
together with a differential $\delta:V_i\to V_{i-1}$, such that
$\delta^2=0$. A map of chain complexes is called a weak equivalence if it
induces isomorphisms in homology. We denote by $V^\sharp$ the graded vector
space underlying $V$, with vanishing differential.

The suspension $\Sigma V_\bull$ of a chain complex $V_\bull$ has components
$(\Sigma V)_n=V_{n-1}$, and differential equal to minus that of
$V_\bull$. By $\Sigma^nV_\bull$, $n\in\Z$, we denote the $n$-fold iterated
suspension of $V_\bull$.

If $V$ is a chain complex and $G$ is a finite group, we denote by $V_G$ the
chain complex of $G$-coinvariants of $V$:
$$
V_G = V/\{gv-v\mid v\in V,g\in G\} .
$$

By $V^*_\bull$, we denote the linear dual of $V_\bull$, with
$V^*_i=(V_{-i})^*$ and with differential $\delta^*:V^*_i\to V^*_{i-1}$ the
adjoint of $\delta:V_{-i+1}\to V_{-i}$. All chain complexes which we consider
in this article have finite dimensional homology.


Denote by $\SS_n$ the group $\Aut\{1,\dots,n\}$ and by $\SS_{n\+}$ the
group $\Aut\{0,1,\dots,n\}$. (This was denoted by $\SS_{n+1}$ in
\cite{cyclic}; we have changed the notation in order to distinguish between
the (isomorphic) groups $\SS_{n\+}=\Aut\{0,\dots,n\}$ and
$\SS_{n+1}=\Aut\{1,\dots,n+1\}$.) An $\SS$-module is a sequence of chain
complexes $\v=\{\v(n)\mid n\ge0\}$, together with an action of $\SS_n$ on
$\v(n)$ for each $n$.

A map of $\SS$-modules is called a weak equivalence if it is a weak
equivalence for each $n$.

\subsection{Operads}
An operad is an $\SS$-module $\pp$ together with bilinear operations
$$
\circ_i : \pp(m) \o \pp(n) \to \pp(m+n-1) , 1\le i\le m ,
$$
satisfying the following axioms.

\subsubsection{} \label{axiom1}
If $\pi\in\SS_m$, $\rho\in\SS_n$, $a\in\pp(m)$ and $b\in\pp(n)$, then
$$
(\pi a)\circ_{\pi(i)}(\rho b) = (\pi\circ_i\rho)(a\circ_ib) ,
$$
where $\pi\circ_i\rho\in\SS_{m+n-1}$ is defined as follows: it permutes the
interval $\{i,\dots,i+n-1\}$ according to the permutation $\rho$, and then
reorders the $m$ intervals
$$
\{1\},\dots,\{i-1\},\{i,\dots,i+n-1\},\{i+n\},\dots,\{m+n-1\} ,
$$
which partition $\{1,\dots,m+n-1\}$, according to the permutation
$\pi$. Explicitly,
$$
(\pi\circ_i\rho)(j) = \begin{cases}
\pi(j) , & \text{$j<i$ and $\pi(j)<\pi(i)$} , \\
\pi(j) + n-1 , & \text{$j<i$ and $\pi(j)>\pi(i)$} , \\
\pi(j-n+1) , & \text{$j\ge i+n$ and $\pi(j)<\pi(i)$} , \\
\pi(j-n+1) + n-1 , & \text{$j\ge i+n$ and $\pi(j)>\pi(i)$} , \\
\pi(i) + \rho(j-i+1) - 1 , & i\le j<i+n .
\end{cases}$$

\subsubsection{} \label{axiom2}
For $a\in\pp(k)$, $b\in\pp(l)$ and $c\in\pp(m)$, and $1\le i<j\le k$,
$$
(a\circ_i b)\circ_{j+l-1}c = (a\circ_j c)\circ_i b .
$$

\subsubsection{} \label{axiom3}
For $a\in\pp(k)$, $b\in\pp(l)$ and $c\in\pp(m)$, and $1\le i\le k$,
$1\le j\le l$,
$$
(a\circ_i b)\circ_{i+j-1}c = a\circ_i(b\circ_j c) .
$$

\subsection{Operads and trees}
We think of an element of $\pp(n)$ as corresponding to a rooted tree with
one vertex, $n$ inputs numbered from $1$ up to $n$ (and one output).
$$\begin{picture}(125,135)(55,665)
\put(120,740){\circle*{4}}
\put(120,800){\line( 0,-1){ 60}}
\put(120,740){\line(-1,-1){ 60}}
\put(120,740){\line( 1,-1){ 60}}
\put(120,740){\line(-2,-3){ 40}}
\put(120,740){\line( 2,-3){ 40}}
\put( 55,665){$1$}
\put( 75,665){$2$}
\put(143,665){$n-1$}
\put(180,665){$n$}
\put(114,685){$\dots$}
\end{picture}$$
The compositions correspond to grafting two such trees together along the
input of the first tree numbered $i$. Axiom \ref{axiom1} expresses the
equivariance of this construction.
$$\begin{picture}(126,195)(125,515)
\put(190,590){\circle*{4}}
\put(180,595){$b$}
\put(190,650){\line( 0,-1){ 60}}
\put(190,590){\line(-1,-1){ 60}}
\put(190,590){\line( 1,-1){ 60}}
\put(190,590){\line( 2,-3){ 40}}
\put(190,590){\line(-2,-3){ 40}}
\put(125,515){$1$}
\put(145,515){$2$}
\put(213,515){$n-1$}
\put(250,515){$n$}
\put(183,535){$\dots$}
\put(190,650){\circle*{4}}
\put(180,655){$a$}
\put(190,650){\line( 0, 1){ 60}}
\put(190,650){\line(-1,-1){ 60}}
\put(190,650){\line( 1,-1){ 60}}
\put(125,575){$1$}
\put(245,575){$m$}
\put(150,590){$\dots$}
\put(215,590){$\dots$}
\put(188,565){$i$}
\end{picture}$$
Axioms \ref{axiom2} and \ref{axiom3} mean that we can construct unambiguous
compositions corresponding to the following two trees respectively:
$$\begin{picture}(160,215)(160,585)
\put(240,800){\line( 0,-1){ 60}}
\put(240,740){\circle*{4}}
\put(230,745){$a$}
\put(240,740){\line(-2,-5){ 31.724}}
\put(240,740){\line( 2,-5){ 31.724}}
\put(240,740){\line( 1,-1){ 80}}
\put(240,740){\line(-1,-1){ 80}}
\put(175,660){$\dots$}
\put(232,660){$\dots$}
\put(290,660){$\dots$}
\put(157,645){$1$}
\put(317,645){$k$}
\put(208,660){\circle*{4}}
\put(198,665){$b$}
\put(208,660){\line(-1,-2){ 30}}
\put(208,660){\line( 1,-3){ 20}}
\put(195,600){$\dots$}
\put(205,635){$i$}
\put(175,585){$1$}
\put(227,585){$l$}
\put(272,660){\circle*{4}}
\put(262,665){$c$}
\put(271,635){$j$}
\put(272,660){\line(-1,-3){ 20}}
\put(272,660){\line( 1,-2){ 30}}
\put(248,585){$1$}
\put(300,585){$m$}
\put(267,600){$\dots$}
\end{picture}
\quad\quad\quad
\begin{picture}(160,235)(160,585)
\put(240,800){\line( 0,-1){ 60}}
\put(240,740){\circle*{4}}
\put(230,745){$a$}
\put(240,740){\line( 1,-1){ 80}}
\put(240,740){\line(-1,-1){ 80}}
\put(240,740){\line( 0,-1){ 60}}
\put(240,680){\circle*{4}}
\put(230,685){$b$}
\put(240,680){\line( 1,-1){ 60}}
\put(240,680){\line( 1,-6){ 08}}
\put(240,680){\line(-1,-1){ 60}}
\put(240,680){\line(-2,-3){ 40}}
\put(185,660){$\dots$}
\put(280,660){$\dots$}
\put(157,645){$1$}
\put(234,650){$i$}
\put(319,645){$k$}
\put(248,632){\circle*{4}}
\put(238,637){$c$}
\put(248,632){\line(-1,-2){ 20}}
\put(248,632){\line( 1,-2){ 20}}
\put(176,605){$1$}
\put(197,605){$2$}
\put(245,610){$j$}
\put(301,605){$l$}
\put(213,620){$\dots$}
\put(268,620){$\dots$}
\put(243,592){$\dots$}
\put(223,577){$1$}
\put(264,577){$m$}
\end{picture}$$
In fact, the axioms imply that the products $\circ_i$ give rise to an
unambiguous definition of composition for any rooted tree \cite{GJ},
\cite{GK}. This point of view will be explained in greater detail, in the
context of modular operads, in Section 2.

\subsection{Cyclic $\SS$-modules}
A cyclic $\SS$-module $\v$ is a sequence of vector space $\v(n)$, with
action of $\SS_{n\+}$ on $\v(n)$. In particular, each vector space $\v(n)$
is a module over the symmetric group $\SS_n$, and over the cyclic group
$\Z_{n\+}$ generated by $\tau_n=(0 1 \dots n)$.

If $\v$ is a cyclic $\SS$-module, and $I$ is a $(k+1)$-element set,
define
$$
\v\(I\) = \Big(
\bigoplus_{\substack{\text{bijections}\\f:\{0,\dots,k\}\to I}} \v(k)
\Bigr)_{\SS_{k\+}} .
$$
This makes $\v$ into a functor from the category of nonempty finite sets
and their bijections into the category of vector spaces. In the case when
$k=n-1$ and $I=\{1,\dots,n\}$, we write $\v\(n\)$ instead of $\v\(I\)$.
Note that $\v\(n\)=\v(n-1)$.

\subsection{Cyclic operads}
If $\pp$ is a cyclic $\SS$-module and $a\in\pp(n)$, let $a^*\in\pp(n)$ be
the result of applying the cycle $(01...n)\in\SS_{n\+}$ to $a$. (Thus, if
$n=1$, this operation exchanges the input and output of $a$; for $n>1$, it
generalizes this.) A cyclic operad \cite{cyclic} is a cyclic $\SS$-module
$\pp$ whose underlying $\SS$-module has the structure of an operad, such
that
\begin{equation}\label{cyclic} 
(a\circ_m b)^* = b^*\circ_1 a^* .
\end{equation}
for any $a\in\pp(m)$, $b\in\pp(n)$. This formula shows that cyclic operads
are a generalization of associative $\ast$-algebras, which are the special
case in which $\pp(n)=0$ for $n\ne1$.

(Cyclic) $\SS$-modules may be defined, in exactly the same way, in any
symmetric monoidal category. The most important case for us will be the
category of chain complexes, whose operads will be called differential
graded operads (abbreviated to dg-operad). Other examples are the category
of topological spaces, giving rise to topological $\SS$-modules and
operads, and the opposite category to the category of chain complexes,
whose operads are called dg-cooperads.

In the remainder of this paper, unless otherwise specified, by an
$\SS$-module, operad or cooperad, we mean a dg $\SS$-module, operad or
cooperad. A map of operads is called a weak equivalence if it is a weak
equivalence of the underlying $\SS$-module.

\subsection{Endomorphism operads and cyclic algebras}
\label{endomorphism-cyclic}
Let $V$ be a chain complex such that its homogeneous subspaces $V_i$ are
finite-dimensional for all $k$. An inner product on $V$ is a non-degenerate
bilinear form $B(x,y)$ such that $B(\delta x,y)+(-1)^{|x|}B(x,\delta y)=0$,
where $\delta$ is the differential of $V$. Such a bilinear form is
symmetric (respectively antisymmetric) if $B(y,x)=(-1)^{|x|\,|y|}B(x,y)$
(resp.\ $B(y,x)=-(-1)^{|x|\,|y|}B(x,y)$), and has degree $k$ if $B(x,y)=0$
unless $|x|+|y|=k$.

Let $V$ be a chain complex with symmetric inner product $B(x,y)$ of degree
$0$. We define a cyclic $\SS$-module $\e[V]$ by putting
$\e[V](n)=V^{\o(n\+)}$, with the natural action of $\SS_{n\+}$. This may
be given the structure of a cyclic operad: if $a\in V^{\o(m\+)}$ and $b\in
V^{\o(n\+)}$, the product $a\circ_i b\in V^{\o(m+n)}$ is defined by
contracting $a\o b$ with the bilinear form $B$, applied to the $i$-th
factor of $a$ and the $0$-th factor of $b$. Using the isomorphism
$V^{\o(n+1)}\cong\Hom(V^{\o n},V)$, the operad underlying this cyclic
operad may be identified with the endomorphism operad of \cite{May} and
\cite{GK}.

A cyclic algebra over a cyclic operad $\pp$ is a chain complex $A$ with
inner product $B$, together with a morphism of cyclic operads
$\pp\to\e[A]$.

\subsection{Examples}

\subsubsection{Stable curves of genus $0$}
Define a topological cyclic operad $\CM_0$ by letting $\CM_0(n)$ be the
moduli space $\CM_{0,n\+}$ of stable curves of genus $0$ with embedding of
$\{0,\dots,n\}$ \cite{Knudsen} (see also \cite{GK}). By definition, a point
of $\CM_{g,n}$ is a system $(C,x_1,\dots,x_n)$, where $C$ is a projective
curve of arithmetic genus $0$, with possibly nodal singularities, $x_i$ are
distinct smooth points, and $C$ has no infinitesimal automorphisms
preserving the points $x_i$ (this amounts to saying that each component of
$C$ minus its singularities and marked points has negative Euler
characteristic). The $\SS_n$-action on $\CM_{0,n}$ is given by renumbering
the punctures. The composition $\circ_i$ takes two pointed curves
$(C,x_0,\dots,x_m)$ and $(D,y_0,\dots,y_n)$ into
$$
\bigl( C\Coprod D/(x_i\sim y_0) ,
x_0,\dots,x_{i-1},y_1,\dots,y_n,x_{i+1},\dots,x_m\bigl) .
$$

\subsubsection{Spheres with holes} \label{spheres}
Define a topological cyclic operad $\Mhat_0$ by letting $\Mhat_0(n)$ be the
moduli space of data $(C,f_0,\dots,f_n)$, where $C$ is a complex manifold
isomorphic to $\CP^1$, and $f_i$ are biholomorphic maps of the unit disk
$$
\Delta = \{z\in \C \mid |z|\le1 \}
$$
into $C$ with disjoint images. The composition $\circ_i$ takes
$(C,f_0,\dots,f_m)$ and $(D,g_0,\dots,g_n)$ into
$$
\Bigl( \bigl( C\setminus
f_i\bigl[\overset{\scriptscriptstyle\circ}\Delta\bigr] \bigr)
\coprod\nolimits_{f_i(t)\sim g_0(t^{-1}),t\in\p\Delta}
\bigl( D\setminus g_0\bigl[\overset{\scriptscriptstyle\circ}\Delta\bigr]
\bigr) ,
f_0,\dots,f_{i-1},g_1,\dots,g_n,f_{i+1},\dots,f_m \Bigr)
$$

Note that by applying the total homology functor $H_\bull(-,\k)$ to the
topological operads $\Mbar_0$ and $\Mhat_0$, we obtain cyclic operads in
the category of graded vector spaces.

\subsubsection{Commutative operad} \label{comm}
This operad, denoted $\Com$, has $\Com\(n\)=\k$ (the trivial representation
of $\SS_n$) for all $n\ge3$, with the obvious composition maps. Cyclic
algebras over $\Com$ are commutative algebras (possibly without unit) with
an invariant scalar product in the ordinary sense: $B(xy,z) = B(x,yz)$.

\subsubsection{Associative operad} \label{ass}
This operad, denoted $\Ass$, has $\Ass\(n\)=\Ind_{\Z_n}^{\SS_n}(\k)$, where
$\k$ is the trivial representation of the cyclic group $\Z_n$. Thus, there
is a natural basis for $\Ass\(n\)$ labelled by the \emph{cyclic orders} of
the set $\{1,\dots,n\}$, that is, the free $\Z_n$-actions on this set. Note
that cyclic orders on $\{0,1,\dots,n\}$ are in bijection with permutations
of $\{1,\dots,n\}$; thus $\Ass\(n\+\)=\Ass(n)$ is free as an
$\SS_n$-module, with generating vector $e_n$ and basis
$\{\sigma\*e_n\mid\sigma\in\SS_n\}$. The composition is determined by the
formulas $e_m\circ_1e_n=e_{m+n-1}$ together with \ref{axiom1}.

A cyclic $\Ass$-algebra is the same as an associative algebra $A$ with an
invariant scalar product. The basis element $\sigma\*e_n\in\Ass(n)$ acts on
$A$ as an $n$-ary operation $(x_1,\dots,x_n)\mapsto x_{\sigma^{-1}(1)}\dots
x_{\sigma^{-1}(n)}$. See \cite{GK}, \cite{cyclic} for more details.

\subsubsection{Lie operad} \label{lie}
This operad, denoted $\Lie$, is determined by the requirement that its
cyclic algebras are Lie algebras with invariant scalar product. Thus
$\Lie(n)=\Lie\(n+1\)$ can be identified, as a module over
$\SS_n\subset\SS_{n\+}$, with the subspace in the free Lie algebra on
generators $x_1,\dots,x_n$ spanned by all Lie monomials containing each
$x_i$ exactly once. The $\SS_n$-module $\Lie(n)$ is isomorphic to the
induced representation $\Ind_{\Z_n}^{\SS_n}(\chi)$, where $\chi$ is a
primitive character of $\Z_n$ (one which takes the generator of $\Z_n$ into
a primitive root of 1).

\section{Modular operads} \label{Graphs}

\subsection{Stable $\SS$-modules}
A stable $\SS$-module is a collection of chain complexes
$$
\{ \v\(g,n\) \mid n,g\ge0 \}
$$
with an action of $\SS_n$ on $\v\(g,n\)$, such that $\v\(g,n\)=0$ if
$2g+n-2\le0$.

A morphism $\v\to\w$ of stable $\SS$-modules is a collection of
equivariant maps of chain complexes $\v\(g,n\)\to\w\(g,n\)$.

We have borrowed the term ``stable'' from the theory of moduli spaces of
curves, since the condition of stability is the same in the two settings.

Any cyclic $\SS$-module $\v$ such that $\v\(n\)=0$ for $n\ge2$ may be
regarded as a stable $\SS$-module by setting:
\begin{equation} \label{stab}
\v\(g,n\) = \begin{cases} \v\(n\) , & g=0 , \\ 0 , & g>0 . \end{cases}
\end{equation}

In the other direction, we have the forgetful functor, which we denote by
$\Cyc$. If $\v$ is a stable $\SS$-module, then $\Cyc(\v)$ is a cyclic
$\SS$-module, and
\begin{equation} \label{Cyc}
\Cyc(\v)\(n\) = \v\(0,n\) .
\end{equation}

A stable $\SS$-module $\v$ has a natural extension to all finite sets $I$:
\begin{equation} \label{(g,I)}
\v\(g,I\) = \Bigl(
\bigoplus_{\substack{\text{bijections}\\f:\{1,\dots,n\}\to I}}
\v\(g,n\) \Bigr)_{\SS_n}
\end{equation}

\subsection{Graphs}
A graph $G$ is a finite set $\Flag(G)$ (whose elements are called flags)
together with an involution $\sigma$ and a partition $\lambda$. (By a
partition of a set, we mean a disjoint decomposition into several
unordered, possibly empty, subsets, called blocks.) We say that two flags
$a,b\in\Flag(G)$ meet if they are equivalent under the partition $\lambda$.

The vertices of $G$ are the blocks of the partition $\lambda$, and the set
of them is denoted $\VV(G)$. The subset of $\Flag(G)$ corresponding to a
vertex $v$ is denoted $\Leg(v)$. Its cardinality is called the valence of
$v$, and denoted $n(v)$.

The edges of $G$ are the pairs of flags forming a two-cycle of $\sigma$,
and the set of them is denoted $\Edge(G)$. The legs of $G$ are the
fixed-points of $\sigma$, and the set of them is denoted $\Leg(G)$. Since a
flag forms either a leg or half an edge, we see that
\begin{equation} \label{1}
\sum_{v\in\VV(G)} n(v) = 2|\Edge(G)| + n .
\end{equation}

\subsection{The geometric realization of a graph} We may associate to a
graph the finite one-dimensional cell complex $|G|$, obtained by taking one
copy of $[0,\half]$ for each flag, and imposing the following equivalence
relation: the points $0\in[0,\half]$ are identified for all flags in a
block of the partition $\lambda$, and the points $\half\in[0,\half]$ are
identified for pairs of flags exchanged by the involution $\sigma$. Thus,
two flags in $G$ meet if and only if their corresponding loci in $|G|$
intersect in a point.

For example, the following corresponds to the set of
flags $\{1,\dots,9\}$, the involution $\sigma=(4 6)(5 7)$ and the
partition $\{1,2,3,4,5\}\cup\{6,7,8,9\}$.
$$\begin{picture}(170,64)(140,645)
\put(200,680){\circle*{4}}
\put(200,680){\line(-3, 1){ 60}}
\put(200,680){\line(-1, 0){ 60}}
\put(200,680){\line(-3,-1){ 60}}
\put(160,700){$\scriptstyle1$}
\put(142,685){$\scriptstyle2$}
\put(160,655){$\scriptstyle3$}
\put(220,680){\circle{40}}
\put(205,700){$\scriptstyle4$}
\put(205,655){$\scriptstyle5$}
\put(230,700){$\scriptstyle6$}
\put(230,655){$\scriptstyle7$}
\put(240,680){\circle*{4}}
\put(240,680){\line( 3, 1){ 60}}
\put(240,680){\line( 3,-1){ 60}}
\put(277,700){$\scriptstyle8$}
\put(277,655){$\scriptstyle9$}
\end{picture}$$

Let $H_i(G)$ be the $i$-th homology group of the cell complex $|G|$ with
coefficients in $\k$, and let $b_i(G)$ be the dimension of $H_i(G)$. Then
$b_0(G)$ is the number of components of $G$, and $b_1(G)$ is the number of
circuits of $G$. The graph $G$ is connected if $b_0(G)=1$.

\subsection{Stable graphs}
A labelled graph is a connected graph $G$ together with a map $g$ from
$\VV(G)$ into the non-negative integers. The value $g(v)$ of this map at a
given vertex $v$ is called the genus of $v$. The genus $g(G)$ of a labelled
graph $G$ is defined by the formula
\begin{equation} \label{g(G)}
g(G) = \sum_{v\in\VV(G)} g(v) + b_1(G) = \sum_{v\in\VV(G)} \bigl(g(v)-1) +
|\Edge(G)| + 1 .
\end{equation}
Adding twice \ref{g(G)} to \ref{1}, we see that
\begin{equation} \label{Euler}
2(g-1) + n = \sum_{v\in\VV(G)} \bigl( 2(g(v)-1) + n(v) \bigr) .
\end{equation}
Likewise, adding three times \ref{g(G)} to \ref{1}, we see that
\begin{equation} \label{3(g-1)}
3(g-1) + n = |\Edge(G)| + \sum_{v\in\VV(G)} \bigl( 3(g(v)-1)+n(v) \bigr) .
\end{equation}
Both of these formulas will be needed later.

A forest is a (labelled) graph of genus $0$; a tree is a connected forest.
This definition is slightly different from the definition of trees in
\cite{cyclic}: here, we do not admit the tree with two legs and no
vertices.

A connected labelled graph is called stable if $2(g(v)-1)+n(v)>0$ for each
vertex $v$. 

If $\v$ is a stable $\SS$-module and $G$ be a stable graph, let $\v\(G\)$
be the tensor product
\begin{equation} \label{v(G)}
\v\(G\) = \bigotimes_{v\in\VV(G)} \v\(g(v),v\) .
\end{equation}

\subsection{Morphisms of graphs}
Let $G_0$ and $G_1$ be two graphs. A morphism $f:G_0\to G_1$ is an
injection $f^*:\Flag(G_1)\to\Flag(G_0)$ such that
\begin{enumerate}
\item $\sigma_0\circ f^*=f^*\circ\sigma_1$, where $\sigma_i$, $i=0,1$, are
the involutions of $\Flag(G_i)$;
\item $\sigma_0$ acts freely on the complement of the image of $f^*$ in
$\Flag(G_0)$ (i.e.\ $G_1$ is obtained from $G_0$ by contracting a subset of
its edges);
\item two flags $a$ and $b$ in $G_1$ meet if and only if there is a chain
$(x_0,\dots,x_k)$ of flags in $G_0$ such that $f^*a=x_0$, $\sigma_0
x_{i-1}$ and $x_i$ meet for all $1\le i\le k$, and $f^*b=\sigma_0 x_k$.
\end{enumerate}
This definition is equivalent to that of Kontsevich-Manin \cite{KM}.

A morphism $f:G_0\to G_1$ defines a surjective cellular map
$|f|:|G_0|\to|G_1|$ which is bijective on the legs.

The preimage of a vertex $v\in\VV(G_1)$ under a morphism $f$, denoted
$f^{-1}(v)$, is the graph consisting of those flags in $G_0$ which are
connected to a flag in $\Leg(v)$ by a chain of edges in $G_0$ contracted by
the morphism. Note that $\Leg(f^{-1}(v))=\Leg(v)$.

A morphism $f:G_0\to G_1$ of labelled graphs is a morphism of the
underlying graphs such that the genus of a vertex $v$ of $G_1$ is equal to
the genus of its inverse image $f^{-1}(v)$ in $G_0$.

Let $\Gamma$ be the category of all stable graphs and their morphisms.

\subsection{Contractions of graphs} \label{catgraph}
Let $G$ be a stable graph and let $I\subset\Edge(G)$ is a subset of its
edges. Then there is a unique stable graph $G/I$ with the following
properties.
\begin{enumerate}
\item $\Flag(G/I)$ is obtained from $\Flag(G)$ by deleting the flags
constituting the edges in $I$;
\item the inclusion $\Flag(G/I)\hookrightarrow\Flag(G)$ is a morphism of
graphs $\pi_{G,I}:G\to G/I$.
\end{enumerate}
The graph $G/I$ is called the contraction of $G$ along the set of edges
$I$. Any morphism $f:G\to G'$ of stable graphs is isomorphic to a morphism
of this form. Note that the realization $|G/I|$ is obtained from $|G|$ by
contracting each edge of $I$ to a point.

If $I$ is a set with just one edge $e$, we will abbreviate $G/\{e\}$ and
$\pi_{G,\{e\}}$ to $G/e$ and $\pi_{G,e}$.

\subsection{The category $\Gamma\(g,n\)$}
Let $\Gamma\(g,n\)$ be the category whose objects are pairs $(G,\rho)$
where $G$ is a stable graph of genus $g$ and $\rho$ is a bijection between
$\Leg(G)$ and the set $\{1,\dots,n\}$, and whose morphisms are morphisms of
stable graphs preserving the labelling $\rho$ of the legs. This category
has a terminal object $\ast_{g,n}$, the graph with no edges and one vertex
$v$ of genus $g$ and valence $n$.

\begin{lemma}
The category $\Gamma\(g,n\)$ has a finite number of isomorphism classes of
objects, which we will denote by $[\Gamma\(g,n\)]$.
\end{lemma}
\begin{proof}
Since the graph $G$ is stable, the integer $3(g(v)-1)+n(v)$ is non-negative
at each vertex. By \ref{3(g-1)}, this gives a bound of $3(g-1)+n$ for the
number of edges of $G$, and hence a bound of $6(g-1)+2n$ for the number of
flags. Since there are a finite number of stable graphs of genus $g$ with a
given number of flags, the lemma follows.
\end{proof}

Denote by $\Aut(G)$ the automorphism group of a graph $G$ in
$\Gamma\(g,n\)$. Observe that trees have no non-trivial automorphisms,
since each vertex is uniquely determined by its distance from each of the
legs (or even any two).

\subsection{The triple of stable graphs}
If $\Cat$ is a category, let $\Iso\Cat$ be the subcategory of isomorphisms
of $\Cat$. We define an endofunctor $\MM$ of the category of stable
$\SS$-modules by the formula
\begin{equation}
\MM\v\(g,n\) = \colim_{G\in\Iso\Gamma\(g,n\)} \v\(G\) \cong
\bigoplus_{G\in[\Gamma\(g,n\)]} \v\(G\)_{\Aut(G)} .
\end{equation}

The functor $\MM$ is a triple; that is, there are natural transformations
$\mu:\MM\MM\v\to\MM\v$ and $\eta:\v\to\MM\v$ making it into a monoid in the
monoidal category of endofunctors of the category of stable
$\SS$-modules. We will now construct these and verify the axioms of a
triple.

We may associate to any category $\Cat$ a simplicial category
$\Iso_\bull\Cat$; the objects of $\Iso_k\Cat$ are diagrams
$$
(f_1,\dots,f_k) = [G_0 \xrightarrow{f_1} G_1 \xrightarrow{f_2} \dots
\xrightarrow{f_{k-1}} G_{k-1} \xrightarrow{f_k} G_k]
$$
in $\Cat$, while the morphisms are isomorphisms of such diagrams. The face
maps are given by the usual formulas
$$
\p_i(f_1,\dots,f_k) = \begin{cases} (f_2,\dots,f_k) , & i=0 , \\
(f_1,\dots,f_{i+1}\*f_i,\dots,f_k) , &
1\le i\le k-1 , \\
(f_1,\dots,f_{k-1}) , & i=k ,
\end{cases}$$
as are the degeneracies,
$$
\sigma_i(f_1,\dots,f_k) = (f_1,\dots,f_i,\Id_{G_i},f_{i+1},\dots,f_k) ,
\quad 0\le i\le k .
$$
In particular, $\Iso_0\Cat=\Iso\Cat$.

Every object of $\Iso_k\Gamma\(g,n\)$ is isomorphic to an object made up of
a sequence of contractions $[G \to G/I_1 \to \dots \to G/I_k]$, where
$G\in\Ob\Gamma\(g,n\)$ is a stable graph and $I_1\subset\dots\subset
I_k\subset\Edge(G)$ is a chain of subsets of $\Edge(G)$.

The proof that $\MM$ is a triple rests on the identity
$$
(\MM^{k+1}\v)\(g,n\) =
\colim_{\substack{[G_0\xrightarrow{f_1}\dots\xrightarrow{f_k}G_k] \\
\in\Iso_k\Gamma\(g,n\)}} \v\(G_0\) .
$$
Equivalently,
\begin{equation} \label{2.?}
(\MM^{k+1}\v)\(g,n\) \cong \colim_{G\in\Iso\Gamma\(g,n\)}
\bigoplus_{I_1\subset\dots\subset I_k\subset\Edge(G)} \v\(G\) .
\end{equation}

The multiplication $\mu:\MM\MM\v\to\MM\v$ of $\MM$ is induced by
$\p_1:\Iso_1\Gamma\(g,n\)\to\Iso_0\Gamma\(g,n\)$, which maps the
contraction $G\to G/I$ to $G$. The unit $\eta:\v\to\MM\v$ of $\MM$ is the
inclusion of the summand $\v\(\ast_{g,n}\)\cong\v\(g,n\)$ of $\MM\v\(g,n\)$
associated to the graph $\ast_{g,n}$ with no edges.

The natural transformations $\mu\MM$ and
$\MM\mu:(\MM^3\v)\(g,n\)\to(\MM^2\v)\(g,n\)$ are induced by the functors
$\p_1,\p_2:\Iso_1\Gamma\(g,n\)\to\Iso_0\Gamma\(g,n\)$, which send the
sequence of contractions $G\to G/I_1\to G/I_2$ respectively to the
contractions $G\to G/I_1$ and $G\to G/I_2$. Since $\p_1\*\p_1=\p_1\*\p_2$,
composing either of these with $\mu$ gives the same natural transformation,
proving associativity of multiplication in the triple $\MM$. It is easy to
see that $\eta$ is a unit.

\subsection{Modular operads}\label{modular}
A modular operad $\a$ is an algebra over the triple $\MM$ in the category
of stable $\SS$-modules. This means that there is a structure map
$\mu:\MM\a\to\a$ such that $\mu\*(\mu\a)=\mu\*(\MM\mu):\MM\a\to\a$, and
$\mu\*(\eta\a)=\Id_\a:\a\to\a$. For example, for any stable $\SS$-module
$\v$, $\MM\v$ is a modular operad, called the free modular operad generated
by $\v$. Modular operads may be considered in any symmetric monoidal
categories.

\subsection{Coherence for modular operads} \label{coherent}
By a modular pre-operad, we mean a stable $\SS$-module $\a$ together with a
structure map $\mu:\MM\a\to\a$. We now give a criterion which determines
when a modular pre-operad is a modular operad.

If $\a$ is a modular pre-operad and $G\in\Ob\Gamma\(g,n\)$, denote by
$\mu_G:\a\(G\)\to\a\(g,n\)$ the $\SS_n$-equivariant map obtained by
composing the universal map
$$
\a\(G\) \to \MM\a\(g,n\) = \colim_{H\in\Iso\Gamma\(g,n\)} \a\(H\)
$$
with the structure map $\mu:\MM\a\(g,n\)\to\a\(g,n\)$. We call this map
composition along the graph $G$. We may use the technique of \ref{(g,I)} to
define maps $\mu_G:\a\(G\)\to\a\(g(G),\Leg(G)\)$ for any stable graph $G$
(no longer requiring that the legs of $G$ be numbered).

Given a morphism $f:G_0\to G_1$ of stable graphs, define a morphism
$\a\(f\):\a\(G_0\)\to\a\(G_1\)$ to be the composition
\begin{multline} \label{a(f)}
\a\(G_0\) = \bigotimes_{u\in\VV(G_0)} \a\(g(u),\Leg(u)\) \cong
\bigotimes_{v\in\VV(G_1)} \a\(f^{-1}(v)\) \\
\xrightarrow{\bigotimes_v\mu_{f^{-1}(v)}} \bigotimes_{v\in\VV(G_1)}
\a\(g(v),\Leg(v)\) = \a\(G_1\) .
\end{multline}

\begin{proposition}
A modular pre-operad $\a$ is a modular operad if and only if the morphisms
$\a\(f\)$ define a functor on the category of stable graphs, that is, if
$$
\a\(f_1f_0\) = \a\(f_1\)\a\(f_0\)
$$
for any two composable morphisms.
\end{proposition}
\begin{proof}
The associativity of the functor $f\mapsto\a\(f\)$ implies that $\a$ is a
modular operad; $\a$ is an $\MM$-algebra precisely when
$\a\(f_1f_0)=\a\(f_1\)\a\(f_0\)$ for all diagrams of the form
$$
G_0 \xrightarrow{f_0} G_1 \xrightarrow{f_1} \ast_{g,n} .
$$

On the other hand, if $\a$ is a modular operad, then given a composable
pair of morphisms $G_0 \xrightarrow{f_0} G_1 \xrightarrow{f_1} G_2$, we see
that $\a\(f_1\)\a\(f_0\)$ is the composition
\begin{multline} \label{!}
\a\(G_0\) = \bigotimes_{u\in\VV(G_0)} \a\(g(u),\Leg(u)\) \cong
\bigotimes_{v\in\VV(G_1)} \a\(f_1^{-1}(v)\) \\
\begin{split}
& \xrightarrow{\bigotimes_v\mu_{f_1^{-1}(v)}} \bigotimes_{v\in\VV(G_1)}
\a\(g(v),\Leg(v)\) \cong
\bigotimes_{w\in\VV(G_2)} \a\(f_2^{-1}(w)\) \\
& \xrightarrow{\bigotimes_w\mu_{f_2^{-1}(w)}} \bigotimes_{w\in\VV(G_2)}
\a\(g(w),\Leg(w)\) = \a\(G_2\) .
\end{split}
\end{multline}
But $\a$ is a modular operad, so that associativity holds for $\a$ applied
to the diagram
$$
(f_1f_0)^{-1}(w) \to f_1^{-1}(w) \to \ast_{g,n} ,
$$
for all $w\in\VV(G_2)$. This allows us to rewrite \ref{!} as
\begin{multline*}
\a\(G_0\) = \bigotimes_{u\in\VV(G_0)} \a\(g(u),\Leg(u)\) \cong
\bigotimes_{w\in\VV(G_2)} \a\((f_1f_0)^{-1}(w)\) \\
\xrightarrow{\bigotimes_w\mu_{(f_1f_0)^{-1}(w)}} \bigotimes_{w\in\VV(G_2)}
\a\(g(w),\Leg(w)\) = \a\(G_2\) ,
\end{multline*}
which is $\a\(f_1f_0\)$.
\end{proof}

\subsection{Endomorphism operads and modular algebras} \label{endomorphism}
Let $V$ be a chain complex with symmetric inner product $B(x,y)$ of degree
$0$. The endomorphism modular operad $\e[V]$ of $V$ has as its underlying
stable $\SS$-module
$$
\e[V]\(g,n\) = V^{\o n} .
$$
The composition maps of $\e[V]$ are defined as follows: if $G$ is a graph,
the vector space $\e[V]\(G\)$ may be identified with $V^{\o\Flag(G)}$, and
the composition map is obtained by contracting elements of this chain
complex with the multlinear form $B^{\o\Edge(G)}$, which contracts with the
factors of $V^{\o\Flag(V)}$ corresponding to the flags which are paired up
as edges of the graph $G$.

It is easily seen that the cyclic operad underlying $\e[V]$ is the
endomorphism cyclic operad introduced in \ref{endomorphism-cyclic}.

An algebra over a modular operad $\a$ is a chain complex $V$ with inner
product $B$, together with a morphism of modular operads $\a\to\e[V]$.

\section{The structure of modular operads} \label{Structure}

In this section, we show that a modular operad is a cyclic operad with
additional structure (a grading by genus and contractions on pairs of legs)
satisfying certain conditions.

\subsection{Cyclic operads and the triple of trees} \label{TT}
For a cyclic $\SS$-module $\v$ we define a cyclic $\SS$-module $\TT\v$
by summing over trees:
$$
\TT\v\(n\) = \bigoplus_{T\in\Gamma\(0,n\)} \v\(T\) .
$$
(In \cite{cyclic}, this triple was denoted $\TT_+$. Since we have no need
for $\TT_-$ in this paper, we omit the subscript from our notation.) The
following result is Theorem (2.2) of \cite{cyclic}.
\begin{theorem} \label{CYCLIC}
The functor $\TT$ is a triple and a cyclic operad $\pp$ with
$\pp(0)=\pp(1)=0$ is the same as an algebra over $\TT$.
\end{theorem}

In particular, we have free cyclic operads $\TT\v$, where $\v$ is a cyclic
$\SS$-module.  Note that we have the following commutative diagram of
triples:
$$\begin{diagram}
\node{\text{stable $\SS$-modules}} \arrow{e,t}{\MM} \arrow{s,l}{\Cyc}
\node{\text{stable $\SS$-modules}} \arrow{s,l}{\Cyc} \\
\node{\text{cyclic $\SS$-modules}} \arrow{e,t}{\TT} \node{\text{cyclic
$\SS$-modules}}
\end{diagram}$$

\subsection{Graded cyclic operads}
A graded cyclic operad is a cyclic operad $\pp$ such that $\pp\(n\)$ has an
$\SS_n$-invariant decomposition
$$
\pp\(n\) = \bigoplus_{g=0}^\infty \pp\(g,n\)
$$
and if $a\in\pp\(g,m\)$ and $b\in\pp\(h,m\)$, then
$a\circ_ib\in\pp\(g+h,n+m-2\)$. We say that $\pp$ is a stable graded cyclic
operad if $\pp\(g,n\)=0$ for $2(g-1)+n\le0$.

\begin{lemma}
If $\a$ is a modular operad, then the cyclic $\SS$-module
$$
\a^\flat\(n\) = \bigoplus_g \a\(g,n\)
$$
is a stable graded cyclic operad.
\end{lemma}
\begin{proof}
If $\v$ is a stable $\SS$-module, the sub-triple of $\MM\v$ induced by
summing over simply connected graphs alone is isomorphic to the triple
$\TT\v^\flat$. It follows that if $\v$ is an $\MM$-algebra, then $\v^\flat$
is a $\TT$-algebra, that is, a cyclic operad. It is clear that it is stable
and graded.
\end{proof}

\subsection{The contraction maps}
Given a finite set $I$ and distinct elements $i,k\in I$, let $G_{g,I}^{ij}$
be the stable graph with $\Flag(G_{g,I}^{ij})=I$, a single vertex with genus
$g$, and a single edge (a loop) joining the flags $i$ and $j$.
$$
G_{g,I}^{ij} =
\begin{picture}(60,40)(-40,-5)
\put(20,0){\circle{40}}
\put(0,0){\circle*{4}}
\put(0,0){\line( 5, 2){20}}
\put(0,0){\line( 4,-1){20}}
\put(0,0){\line(-2, 1){20}}
\put(0,0){\line(-1, 0){20}}
\put(0,0){\line(-2,-1){20}}
\put(-2,13){$i$}
\put(-2,-20){$j$}
\end{picture}$$
\vskip0.3in
If $\a$ is a modular operad, denote by $\xi_{ij}$ the composition map
$$
\mu_{G_{g,I}^{ij}} : \a\(G_{g,I}^{ij}\) \cong \a\(g,I\) \to
\a\(g+1,I\setminus\{i,j\}\) .
$$
(Here we make use of the notation \ref{(g,I)}.) We call $\xi_{ij}$ the
contraction map. These maps are equivariant, in the sense that for any
bijection $\sigma:I\to J$ of finite sets and $i,j\in I$,
\begin{equation} \label{equivariant}
\xi_{\sigma(i)\sigma(j)}\*\sigma = \sigma\*\xi_{ij} , \quad
\text{on $\in\a\(g,I\)$} .
\end{equation}

We now determine the coherence relations that the contractions $\xi_{ij}$
on a stable graded cyclic operad must satisfy in order for them to define a
modular operad structure.
\begin{theorem}
Let $\a$ be a stable graded cyclic operad with contraction maps
$$
\xi_{ij} : \a\(g,I\) \to \a\(g+1,I\setminus\{i,j\}\) ,
$$
equivariant in the sense of \ref{equivariant}. These data determine a
modular operad if and only if the following coherence conditions are
satisfied:

\begin{description}
\item[1)] For any finite set $I$ and distinct elements $i,j,k,l\in I$,
$$
\xi_{ij}\circ \xi_{kl} = \xi_{kl}\circ \xi_{ij} .
$$
\end{description}
The remaining conditions concern the composition $a\circ_mb$ of $a\in\a(m)$
and $b\in\a(n)$:
\begin{description}
\item[2)] \hfil $\xi_{12}(a\circ_m b) = (\xi_{12}a)\circ_m b$; \hfil
\item[3)] \hfil $\xi_{m,m+1}(a\circ_m b) = a\circ_m(\xi_{12}b)$; \hfil
\item[4)] \hfil $\xi_{m-1,m}(a\circ_m b) = \xi_{m+n-2,m+n-1}(a\circ_{m-1}b^*)$. \hfil
\end{description}
\end{theorem} 
\begin{proof}
``Only if'': Let $\a$ be a modular operad. By \ref{coherent}, we obtain a
functor $f\mapsto\a\(f\)$ from the category $\Gamma$ of stable graphs to
$\Cat$. If $e,e'$ are two edges in a stable graph $G$, the contractions of
$e$ and $e'$ commute in $\Gamma$, in the sense that
$$
\pi_{G/e,e'} \* \pi_{G,e} = \pi_{G/e',e} \* \pi_{G,e'} =
\pi_{G,\{e,e'\}} : G \to G/\{e,e'\} .
$$
We now obtain relations 1-4) in the statement of the theorem by evaluating
the functor $\a\(f\)$ on these identities, for all stable graphs with two
edges. Indeed, a graph with two edges has one of the following forms:
\begin{gather*}
\text{a)} {\epsfxsize=150pt\epsffile{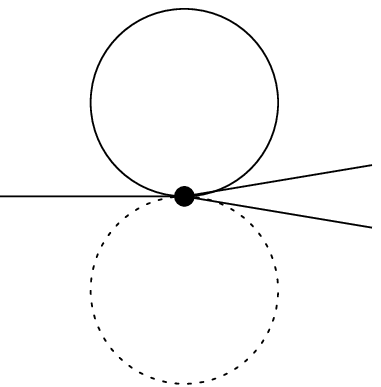}}
\text{b)} {\epsfxsize=150pt\epsffile{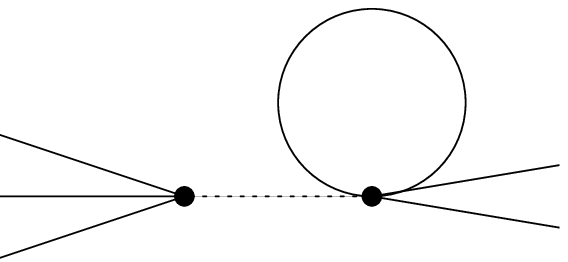}} \\[0.5in]
\text{c)} {\epsfxsize=150pt\epsffile{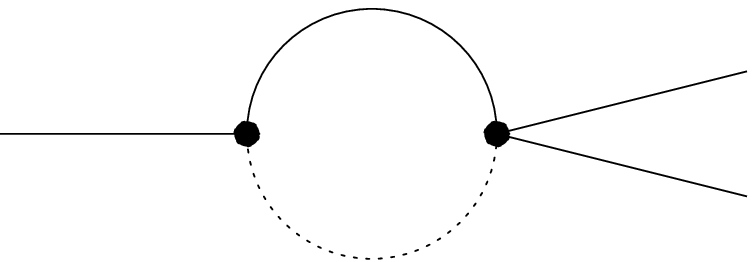}}
\text{d)} {\epsfxsize=150pt\epsffile{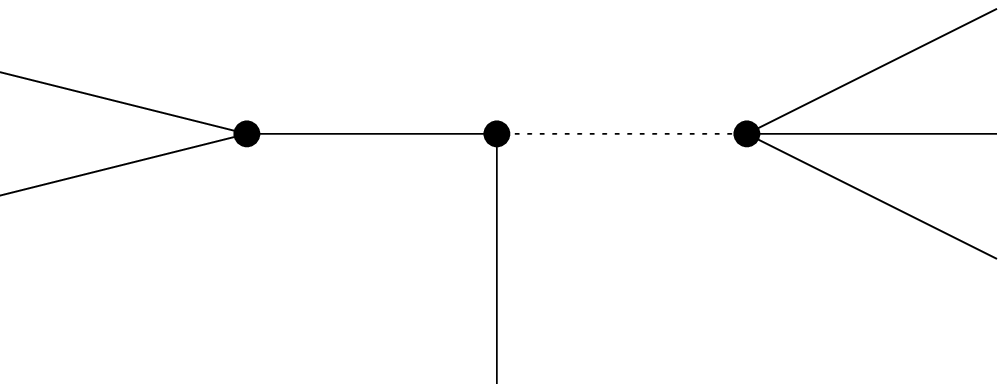}} \\[0.3in]
\end{gather*}
Graphs of type a) give rise to the relations of type 1) in the statement of
the theorem, graphs of type b) to the relations of type 2) and 3) and
graphs of type c) to relations of type 3). (Graphs of type d) of course
imply that $\a^\flat$ is a (graded) cyclic operad.)

``If'': Let $\a$ be a stable graded cyclic operad equipped with contraction
maps $\xi_{ij}$ as in the statement of the theorem. We will construct a
functor, which we still denote by $\a$, from the category $\Gamma$ of
stable graphs to $\Cat$. On objects $G$ of $\Gamma$ (stable graphs), we
define $\a\(G\)$ as in \ref{v(G)}. If $\phi:G_0\to G_1$ is an isomorphism
of stable graphs, we define $\phi:\a\(G_0\)\to\a\(G_1\)$ in the evident
way: this definition is functorial on the subcategory $\Iso\Gamma$.

Next, we turn to the case where $f=\pi_{G,e}:G\to G/e$ is a
contraction. There are two sub-cases:

1) The edge $e$ has two distinct ends, the vertices $v$ and $v'$: Let $H$
be the graph whose flags are the legs of $v$ and $v'$. Then $H$ is a tree
with one edge $e$ and two vertices $v$ and $v'$, and $H/e$ has a single
vertex $\bar{v}$. We have
\begin{align*}
\a\(G\) &\cong \a\(H\) \o \bigotimes_{w\ne v,v'} \a\(g(w),\Leg(w)\) , \\
\a\(G/e\) &\cong \a\(g(\bar{v}),\Leg(\bar{v})\) \o \bigotimes_{w\ne v,v'}
\a\(g(w),\Leg(w)\) ,
\end{align*}
and we define
$$
\a\(\pi_{G,e}\) = \mu_H \o \bigotimes_{w\ne v,v'} \Id_{\a\(g(w),\Leg(w)\)} ,
$$
where $\mu_H$ is the composition along $H$ in the graded cyclic operad
$\a^\flat$.

2) The edge $e$ has one end, the vertex $v$: Let $i,j\in\Leg(v)$ be the two
flags of $e$, and let $\bar{v}$ be the image of the vertex $v$ in $G/e$. We
have
\begin{align*}
\a\(G\) &\cong \a\(g(v),\Leg(v)\) \o \bigotimes_{w\ne v}
\a\(g(w),\Leg(w)\) , \\
\a\(G/e\) &\cong \a\(g(\bar{v})+1,\Leg(\bar{v})\setminus\{i,j\}\) \o
\bigotimes_{w\ne v} \a\(g(w),\Leg(w)\) ,
\end{align*}
and we define
$\a\(\pi_{G,e}\)=\xi_{ij}\o\bigotimes_{w\ne v}\Id_{\a\(g(w),\Leg(w)\)}$.

Now let $f:G_0\to G_1$ be a general morphism of stable graphs, and let $I$
be the set of edges of $G$ contracted by $f$. The morphism $f$ decomposes
as a composition
$$
G_0 \xrightarrow{\pi_{G_0,I}} G_0/I \xrightarrow{\phi} G_1 ,
$$
where $\phi:G_0/I\to G_1$ is an isomorphism. Choosing an ordering
$\{e_1,\dots,e_k\}$ of the edges in $I$, we obtain a factorization
$$
G_0 \xrightarrow{\pi_{G_0,e_1}} G_0/e_1
\xrightarrow{\pi_{G_0\setminus\{e_1\},e_2}} G_0/\{e_1,e_2\}
\xrightarrow{\pi_{G_0\setminus\{e_1,e_2\},e_3}}
\dots \xrightarrow{\pi_{G_0\setminus\{e_1,\dots,e_{k-1}\},e_k}} G_0/I
\xrightarrow{\phi} G_1 ,
$$
where each morphism is a contraction along one edge except the last, which
is an isomorphism. We define
$$
\a\(f\) = \a\(\phi\) \a\(\pi_{G_0\setminus\{e_1,\dots,e_{k-1}\},e_k}\)
\dots \a\(\pi_{G_0,e_1}\) .
$$

We must prove that this definition of $\a\(f\)$ is independent of the
ordering of the elements of $I$. It suffices to prove the product does not change
if we interchange two consecutive edges $e_i$ and $e_{i+1}$. If the two
edges do not meet, this is evident. If they do meet, then they form a
stable graph with two edges, whose topology is one of the four types a-d)
catalogued above. From conditions 1-4) and equivariance \ref{equivariant},
it follows that composition is well-defined along every graph with two
edges, regardless of the numbering of its legs.

Thus, under the hypotheses of the theorem, we have defined a functor
$f\mapsto\a\(f\)$ on the category of stable graphs, using the graded cyclic
operad structure $\a^\flat$ and the contractions $\xi_{ij}$. It remains to
show that these functors are related by \ref{a(f)} to the underlying
modular pre-operad structure $\mu:\MM\a\to\a$ associated to the special
morphisms of graphs $G\to\ast_{g,n}$. This is evident if we order the
vertices of $G_1$, and then order $I$ in a compatible fashion; the
decomposition of $\a\(f\)$ obtained from this ordering clearly corresponds
to \ref{a(f)}.
\end{proof}

\section{Twisted modular operads}

In this section, we introduce twisted triples $\MM_\local$, which will be
used in the next section in the construction of the Feynman transform .

\subsection{Hyperoperads}
A hyperoperad $\local$ in a symmetric monoidal category $\Cat$ is a
collection of functors from the categories $\Iso\Gamma\(g,n\)$, $g,n\ge0$,
of stable graphs and their isomorphisms to $\Cat$, together with the
following data.

\subsubsection{}
To each morphism $f:G_0\to G_1$ of $\Gamma\(g,n\)$ is assigned a morphism
in $\Cat$
$$
\nu_f : \local(G_1) \o \bigotimes_{v\in\VV(G_1)} \local(f^{-1}(v)) \to
\local(G_0) ,
$$
natural with respect to isomorphisms.

\subsubsection{} \label{hyper.b}
If $\ast_{g,n}$ is the graph in $\Gamma\(g,n\)$ with no edges,
$\local(\ast_{g,n})\cong\1$, the unit object of $\Cat$.

These data are required to satisfy the following conditions.

\subsubsection{} \label{hyper.1}
Given a sequence of morphisms $G_0\xrightarrow{f_1}G_1\xrightarrow{f_2}G_2$
in $\Gamma\(g,n\)$, the following diagram commutes:
$$\begin{diagram}
\node{\local(G_2) \o \bigotimes_{v\in\VV(G_2)}
\local\bigl(f_2^{-1}(v)\bigr) \o \bigotimes_{v\in\VV(G_1)}
\local\bigl(f_1^{-1}(v)\bigr)} \arrow{e,t}{\nu_{f_2}}
\arrow{s,l,=}{\cong}
\node{\local(G_1) \o \bigotimes_{v\in\VV(G_1)} \local\bigl(f_1^{-1}(v)\bigr)}
\arrow[2]{s,r}{\nu_{f_1}} \\
\node{\local(G_2) \o \bigotimes_{v\in\VV(G_2)} \Bigl(
\local\bigl(f_2^{-1}(v)\bigr) \o {\textstyle\bigotimes_{w\in f_2^{-1}(v)}}
\local\bigl(f_1^{-1}(w)\bigr) \Bigr)}
\arrow{s,l}{\scriptstyle\local(G_2)\o\bigotimes_{v\in\VV(G_2)}\nu_{f_1|f_2^{-1}(v)}}
\\
\node{\local(G_2) \o \bigotimes_{v\in\VV(G_2)}
\local\bigl((f_2\*f_1)^{-1}(v)\bigr)}
\arrow{e,b}{\nu_{f_2\*f_1}} \node{\local(G_0)}
\end{diagram}$$
Here, $f_1|f_2^{-1}(v)$ denotes the restriction of $f_1$ to the subgraph
$f_2^{-1}(v)$ of $G_1$.

\subsubsection{} \label{hyper.2}
If $f:G_0\to G_1$ is an isomorphism, the following diagram commutes:
$$\begin{diagram}
\node{\textstyle\local(G_1)\o\bigotimes_{v\in\VV(G_1)}\local(f^{-1}(v))}
\arrow{se,=} \arrow{e,t}{\nu_f} \node{\local(G_0)} \arrow{s,r}{\local(f)}\\
\node[2]{\local(G_1)}
\end{diagram}$$

\subsection{Modular $\local$-operads}
If $\local$ is a hyperoperad, define an endofunctor $\MM_\local$ on the
category of stable $\SS$-modules by the formula
$$
\MM_\local\v\(g,n\) = \bigoplus_{G\in\Iso\Gamma\(g,n\)} \local(G) \o
\v\(G\) .
$$
We show that $\MM_\local$ is a triple by imitating the proof that $\MM$
is. The unit of the triple is again defined by the inclusion of the summand
associated to the graph $\ast_{g,n}$ with no edges in $\Gamma\(g,n\)$: we
may identify $\v(g,n\)$ with $\local(\ast_{g,n})\o\v\(\ast_{g,n}\)$, and
$\local(\ast_{g,n})\cong\1$ by \ref{hyper.b}.

We have the identity
$$
(\MM_\local^2\v)\(g,n\) =
\colim_{\substack{[G_0\stackrel{f_1}{\to}G_1] \\
\in\Iso_1\Gamma\(g,n\)}} \Bigl\{ \local(G_1) \o \bigotimes_{v\in\VV(G_1)}
\local(f_1^{-1}(v)) \o \v\(G_0\) \Bigr\} .
$$
Using the hyperoperad structure maps $v_f$, it is easy to define a natural
transformation $\mu$ from $(\MM_\local^2\v)\(g,n\)$ to
$(\MM_\local\v)\(g,n\)$. By \ref{hyper.2}, we see that $\eta$ is a unit for
this multiplication.

We also have the identity
$$
(\MM_\local^3\v)\(g,n\) =
\colim_{\substack{[G_0\stackrel{f_1}\to G_1\stackrel{f_2}\to G_2] \\
\in \Iso_2\Gamma\(g,n\)}} \Bigl\{ \local(G_2) \o \bigotimes_{v\in\VV(G_2)}
\local(f_2^{-1}(v)) \o \bigotimes_{v\in\VV(G_1)} \local(f_1^{-1}(v))
\o \v\(G_0\) \Bigr\} .
$$
By \ref{hyper.1}, we see that the multiplication $\mu$ is associative.

A modular $\local$-operad is an algebra over the triple $\MM_\local$.

\subsection{Cocycles}
If $\Cat$ is a symmetric monoidal category, an object $L$ is said to be
invertible if there is an object $L^{-1}$ and an isomorphism $\1\cong L\o
L^{-1}$.

A cocycle is a hyperoperad $\local$ such that $\local(G)$ is invertible for
all stable graphs $G$, and such that the morphisms $\local_f$ associated to
morphisms of stable graphs $f:G_0\to G_1$ are isomorphisms. The inverse of
a cocycle $\local$ is again a cocycle, which we denote by $\local^{-1}$.

\subsection{Coboundaries} \label{coboundary}
Let $\lam$ be an $\SS$-module such that each object $\lam\(g,n\)$ is
invertible. Tensoring with $\lam$ defines a functor on $\SS$-modules, which
we denote by $\v\mapsto\lam\v$. There is a natural structure of a cocycle
on
$$
\local_\lam(G) = \lam\(g,n\) \o \bigotimes_{v\in\VV(G)}
\lam\(g(v),n(v)\)^{-1} ,
$$
and a natural isomorphism of triples
$\MM_{\local\o\local_\lam}\cong\lam\circ\MM_\local\circ\lam^{-1}$. We call
this cocycle the coboundary of $\lam$. It follows that if $\local$ is a
cocycle, the functor $\lam$ induces an equivalence between the category of
modular $\local$-operads and the category of modular
$\local\o\local_\lam$-operads.

Let $\local_\Susp$ be the coboundary associated to the invertible stable
$\SS$-module $\Susp$ given by
$$
\Susp\(g,n\) = \Sigma^{-2(g-1)-n} \sgn_n ,
$$
where $\sgn_n$ is the alternating character of $\SS_n$. Equation
\ref{Euler} shows that $\local_\Susp$ is concentrated in degree $0$.

The functor $\v\mapsto\Susp\v$ is called suspension. Since
$\local_\Susp^2\cong\1$, the double suspension of a modular operad is a
modular operad. Note that the suspension of cyclic $\SS$-modules,
considered as stable $\SS$-modules, coincides with the definition of
suspension on cyclic $\SS$-modules \cite{cyclic}, given by the formula
$\Lambda\v(n)=\Sigma^{1-n}\sgn_{n+1}\o\v(n)$.

Two further coboundaries which will be of of interest are associated to the
invertible stable $\SS$-modules $\Pp\(g,n\)=\Sigma^{-6(g-1)-2n}\k$ and
$\Sigma\(g,n\)=\Sigma\k$.

\subsection{Determinants}
For a finite-dimensional vector space $V$ of dimension $n$, let $\Det(V)$
be the graded vector space
$$
\Det(V) = \Sigma^{-n} \Lambda^nV ;
$$
this is the one-dimensional top exterior power of $V$, concentrated in
degree $-n$. If $S$ is a finite set, let $\Det(S)=\Det(\k^S)$. Observe that
there is a natural isomorphism
\begin{equation} \label{invert-det}
\Det(S)^2 \cong \Sigma^{-2|S|}\k .
\end{equation}

\begin{lemma} \label{det}
Given a collection of vector spaces $(V_i)_{i\in I}$, there is natural
identification
$$
\Det\left(\oplus_i V_i\right) \simeq \bigotimes_{i\in I} \Det(V_i) .
$$
\end{lemma}

\subsection{The dualizing cocycle} \label{can}
By \ref{det}, we see that $\can(G)=\Det(\Edge(G))$ is a cocycle,
which we call the dualizing cocycle. Given a cocycle $\local$, we denote
the cocycle $\can\o\local^{-1}$ by $\local^\Dual$ and call it the dual of
$\local$. This duality will be important in the definition of the Feynman
transform for modular operads.
\begin{proposition} \label{Pi}
There is a natural isomorphism of cocycles $\can^2\cong\local_\Pp$.
\end{proposition}
\begin{proof}
There is a natural isomorphism $\local_\Pp(G)\cong\Sigma^{-2\ell}\k$, where
$$
\ell = 3(g-1) + n - \sum_{v\in\VV(G)} \bigl( 3(g(v)-1) + n(v) \bigr) .
$$
By \ref{3(g-1)}, we see that $\ell=|\Edge(G)|$, from which the result
follows.
\end{proof}

\subsection{The orientation cocycle} Let $\OR(e)$ be the orientation line
of an edge $e$ in a graph $G$, that is, the determinant
$\Sigma^2\Det(\{s,t\})$, where $s$ and $t$ are the pair of flags making up
the edge $e$. The orientation cocycle $\twist(G)$ of a graph $G$ is the
one-dimensional vector space
$$
\twist(G) = \Det\Bigl( \bigoplus_{e\in\Edge(G)} \OR(e) \Bigr) .
$$
\begin{proposition} \label{Lambda}
There is a natural isomorphism $\can\cong\twist\o\local_\Susp$.
\end{proposition}
\begin{proof}
If $x$ and $y$ are two independent elements of a vector space $V$, denote
by $\Sigma^{-1}x\.\Sigma^{-1}y$ the corresponding element of
$\Det(\operatorname{Span}\{x,y\})$. If $s$ and $t$ are the two flags making
up an edge $e$, then $\OR(e)$ is spanned by
$\Sigma^2(\Sigma^{-1}s\.\Sigma^{-1}t)$, and thus $\Det(\OR(e))$ is spanned
by a vector $\Sigma(\Sigma^{-1}s\.\Sigma^{-1}t)$. We may identify this with
the element $\Sigma e\o(\Sigma^{-1}s\.\Sigma^{-1}t)$ of
$\Det(\{e\})^{-1}\o\Det(\{s,t\})$. Tensoring over all edges of $G$, we
obtain a natural isomorphism
$$
\twist(G) \cong \Det(\Edge(G))^{-1} \o \Det(\Flag(G)) \o \Det(\Leg(G))^{-1} .
$$
(Here, we use the fact that $\Det(\Flag(G))\o\Det(\Leg(G))^{-1}$ is the
$\Det$ of the set of \textit{internal} flags of $G$, those which are not
legs.)

Thus, it remains to shows that
$$
\local_\Susp(G) \cong \Det(\Edge(G))^2 \o \Det(\Flag(G))^{-1} \o
\Det(\Leg(G)) .
$$
Since $\local_\Susp$ is concentrated in degree $0$, we see that
$$
\local_\Susp(G) \cong \Sigma^{-2|\Edge(G)|} \Det(\Flag(G))^{-1} \o
\Det(\Leg(G)) .
$$
The proposition now follows from \ref{invert-det}, which shows that
$\Det(\Edge(G))^2\cong\Sigma^{-2|\Edge(G)|}$.
\end{proof}

Modular $\twist$-operads admit a notion of algebra parallel to that for
modular operads, as is shown by the following proposition.
\begin{proposition}
Let $V$ be a chain complex with antisymmetric inner product $B(x,y)$ of
degree $-1$. Define the stable $\SS$-module $\e[V]$ of endomorphisms:
$$
\e[V]\(g,n\) = V^{\o n} .
$$
There is a modular $\twist$-operad structure on $\e[V]$.
\end{proposition}
\begin{proof}
For a graph $G\in\Gamma\(g,n\)$, the composition
$\e[V]\(G\)\o\twist(G)\to\e[V]\(g,n\)$ is defined in the same way as in the
untwisted case: we identify $\e[V]\(G\)$ with $V^{\o\Flag(G)}$ and contract
with $B^{\o\Edge(G)}$. The resulting map is well-defined, and invariant
under the action of the groups $\Aut(G)$ and $\SS_n$: the antisymmetry of
$B$ is needed since reversing an edge $e\in\Edge(G)$ changes the sign of
$\twist(G)$, while the degree of $B$ must be $-1$, since supressing an edge
$e\in\Edge(G)$ changes the degree of $\Det(G)$ by $1$.
\end{proof}

\subsection{The determinant of a graph} \label{DET}
The determinant of a graph $G$ is defined to be $\Det(G)=\Det(H_1(G))$.
\begin{proposition} \label{Det}
There is a natural isomorphism
$$
\Det \cong \twist \o \local_\Sigma^{-1}
\cong \can \o \local_\Susp^{-1} \o \local_\Sigma^{-1} .
$$
In particular, $\Det$ is a cocycle.
\end{proposition}
\begin{proof}
This follows from applying \ref{det} to the exact sequence of vector spaces
arising from the complex of cellular chains of the graph $G$,
$$
0 \to H_1(G) \to \bigoplus_{e\in\Edge(G)}\OR(e) \to \k\o\VV(G) \to
H_0(G)\cong\k \to 0 .
\qed$$
\def\qed{}
\end{proof}

Since $\Det(H_1(G))$ is trivial when $G$ is a tree, we see that a cyclic
operad may be considered as a modular $\local$-operad for either the
trivial cocycle $\local=\1$ or the determinant cocycle $\local=\Det$.

\section{The Feynman transform of a modular operad}

In this section, we define a functor $\FF_\local$ from the category of dg
modular $\local$-operads to the category of dg modular
$\local^\Dual$-operads, where we recall that
$\local^\Dual=\can\o\local^{-1}$, and $\can$ is the dualizing cocycle
\ref{can}. We call it the Feynman transform, since $\FF_\local\a$ is a sum
over graphs, as is Feynman's expansion for amplitudes in quantum field
theory.

The most important properties of $\FF_\local$ are that it is a homotopy
functor, in the sense that it maps weak equivalences to weak equivalences,
and that it has homotopy inverse $\FF_{\local^\Dual}$: that is, there is a
natural transformation from $\FF_{\local^\Dual}\FF_\local$ to the identity
functor such that for any modular operad $\a$,
$\FF_{\local^\Dual}\FF_\local\a\to\a$ is a weak equivalence. In this way,
we see that the homotopy categories of modular $\local$-operads and modular
$\local^\Dual$-operads are equivalent.

In the special case where $\local=\1$ is the trivial cocycle, we denote
$\FF_\1$ by $\FF$, and $\FF_{\1^\Dual}$ by $\FF^{-1}$.

\subsection{Definition of the Feynman transform} \label{FF}
As a stable $\SS$-module, but ignoring differentials, $\FF_\local\a$ equals
$\MM_{\local^\Dual}\a^*$, the underlying stable $\SS$-module of the free
modular $\local^\Dual$-operad generated by the linear dual $\a^*$ of $\a$.
The differential $\delta_{\FF_\local\a}$ is the sum
$\delta_{\FF_\local\a}=\delta_{\a^*}+\p$, where $\delta_{\a^*}$ is the
differential on $\MM_{\local^\Dual}\a^*$ induced by the differential on
$\a^*$, and $\p$ is defined as follows.

If $G$ is a stable graph and $e$ is an edge of $G$, the adjoint of the
structure map of the morphism $\pi_{G,e}:G\to G/e$ is a map
$$
\bigl(\mu_{\pi_{G,e}}\bigr)^* : \local(G/e)^* \o \a\(G/e\)^* \to
\local(G)^* \o \a\(G\)^*
$$
of degree $0$. There is natural map $\eps_e:\can(G/e)\to\can(G)$, given by
tensoring with the natural basis element $e$ of
$\can(\{e\})=\Det(\{e\})$. Tensoring these two maps together, we obtain a
map
$$\begin{diagram}
\node{\can(G/e) \o \local(G/e)^* \o \a\(G/e\)^*} \arrow{s,=}
\arrow[4]{e,t}{\eps_e\o\bigl(\mu_{\pi_{G,e}}\bigr)^*}
\node[4]{\can(G)\o\local(G)^*\o\a\(G\)^*} \arrow{s,=} \\
\node{\local^\Dual(G/e) \o \a\(G/e\)^*} \arrow[4]{e,t,..}{\delta_{G,e}}
\node[4]{\local^\Dual(G)\o\a\(G\)^*}
\end{diagram}$$
of degree $-1$.

Recall that $\MM_{\local^\Dual}\a\(g,n\)$ is the sum of complexes
$\bigl(\can(G)\o\local(G)^*\o\a\(G\)^*\bigr)_{\Aut(G)}$ over isomorphism
classes of stable graphs $G\in\Ob\Gamma\(g,n\)$. Given two stable graphs
$G$ and $H$, define the matrix element
$$
\bigl( \can(G) \o \a\(G\)^* \bigr)_{\Aut(G)} \xrightarrow{\p_{G,H}}
\bigl( \can(H) \o  \a\(H\)^* \bigr)_{\Aut(H)}
$$
to be the sum of the maps $\delta_{H,e}$ over all edges $e$ of $H$ such
that $G\cong H/e$; in particular, $\p_{G,H}$ vanishes if
$|\Edge(G)|\ne|\Edge(H)|-1$. The term $\delta_{H,e}$ does not depend on the
isomorphism of $H/e$ with $G$ which is used, since any two such
isomorphisms differ by an automorphism of $G$, and we have taken
coinvariants with respect to $\Aut(G)$.
\begin{theorem}
1) The map $\delta_{\FF_\local\a}$ has square zero.

\noindent 2) The pair
$(\FF_\local\a=\MM_{\local^\Dual}\a^*,\delta_{\FF_\local\a})$ is a modular
$\local^\Dual$-operad of chain complexes.

\noindent 3) The Feynman transform $\FF_\local$ is a homotopy functor: if
$f:\a\to\b$ is a weak equivalence of modular $\local$-operads, then so is
$\FF_\local f:\FF_\local\b\to\FF_\local\a$.
\end{theorem}
\begin{proof}
The matrix element
$$
\bigl(\can(G)\o\a\(G\)^*\bigr)_{\Aut(G)} \xrightarrow{(\p^2)_{G,K}}
\bigl(\can(K)\o\a\(K\)^*\bigr)_{\Aut(K)}
$$
of $\p^2$ is a sum over pairs $(e_1,e_2)$ of distinct edges of $K$ such
that $G\cong K/\{e_1,e_2\}$. The exchange $(e_1,e_2)\mapsto(e_2,e_1)$ is a
fixed-point free involution on the set of such pairs. The respective
contributions $\delta_{K/e_1,e_2}\*\delta_{K,e_1}$ and
$\delta_{K/e_2,e_1\}}\*\delta_{K,e_2}$ to $\p^2$ cancel, since the two
isomorphisms
$$
\can(K) \cong \can(G) \o \can(\{e_1,e_2\})
$$
in their definition are negatives of each other, showing that $\p^2=0$.

It is clear that $\p\circ\delta_{\a^*}+\delta_{\a^*}\circ\p=0$, since the
differential in $\a$ which induces $\delta_{\a^*}$ is compatible with the
structure maps of the modular $\local$-operad $\a$. Together, these results
show that $\delta_{\FF_\local\a}$ has square zero, proving part 1.

It is obvious that the internal differential $\delta_{\a^*}$ is compatible
with the modular $\local^\Dual$-operad structure. To prove part 2), it
remains to show that the differential $\p$ is compatible with the structure
maps of the modular $\local$-operad $\a$,
\begin{gather*}
\local^\Dual(G)\o\MM_{\local^\Dual}\a^*\(G\) \cong
\colim_{\substack{[f:G'\to G] \\ \in\Iso\Gamma\(g,n\)/G}} \local^\Dual(G) \o
\bigotimes_{v\in\VV(G)} \bigl( \local^\Dual(f^{-1}(v)) \o \a^*\(f^{-1}(v)\)
\bigr) \\ \xrightarrow{\quad\mu_G\quad} \a\(g,n\) ;
\end{gather*}
recall that $\Gamma\(g,n\)/G$ is the comma category, whose objects are
morphisms $f:G'\to G$ in $\Gamma\(g,n\)$ with target $G$. The differential
induced by $\p$ on $\FF_\local\a\(G\)$ is a sum. index by the vertices $v$
of $G$, of terms which are themselves a sum, over the vertices $u$ of the
graph $f^{-1}(v)$, of all ways of inserting an edge at $u$. This is clearly
the same as summing over all ways of inserting an edge at all the vertices
of $G'$. On application of the structure map $\mu_G$, this goes into the
differential $\p$ of $\FF_\local\a\(g,n\)$, showing that $\p$ is compatible
with the modular $\local^\Dual$-operad structure on $\FF_\local\a$.

Part 3) is proved by considering a spectral sequence associated to the cone
of the map $\FF_\local f$: we filter by number of edges, in such a way that
the $E^1$-term of the spectral sequence equals the cone of the map
$\FF_\local H_*(f)$, which is zero by hypothesis. The convergence of the
spectral sequence follows from the fact that $\FF_\local\a\(g,n\)$ and
$\FF_\local\b\(g,n\)$ have contributions from a finite number of graphs, so
that the spectral sequence is uniformly bounded in one direction.
\end{proof}

\subsection{The homotopy inverse of the Feynman transform}
Under our hypothesis that $\a\(g,n\)$ is finite dimensional in each degree,
there are isomorphisms of stable $\SS$-modules
$$
\FF_{\local^\Dual}\FF_\local\a \cong
\MM_\local\bigl(\MM_{\local^\Dual}\a^*\bigr)^* \cong
\MM_\local\MM_{\can^{-1}\o\local}\a ,
$$
which shows that $\bigl(\FF_{\local^\Dual}\FF_\local\a\bigr)\(g,n\)$ is a
colimit over $[G_0\xrightarrow{f}G_1]\in\Iso\Gamma_1\(g,n\)$ of the functor
$$
[G_0\xrightarrow{f}G_1] \mapsto \local(G_1) \o \bigotimes_{v\in\VV(G_1)}
\bigl( \can(f^{-1}(v))^{-1} \o \local(f^{-1}(v)) \bigr) \o \a\(G_0\) .
$$
The summand associated to the object
$[\ast_{g,n}\xrightarrow{\Id}\ast_{g,n}]$ is isomorphic to $\a\(g,n\)$. Let
$\tau:\FF_{\local^\Dual}\FF_\local\a\to\a$ be the map induced by projection
onto this summand. We now arrive at the main result of this section.

\begin{theorem} \label{Feynman}
If $\a$ is a modular $\local$-operad, the canonical map
$\tau:\FF_{\local^\Dual}\FF_\local\a\to\a$ is a weak equivalence, that is,
induces an isomorphism on homology.
\end{theorem}
\begin{proof}
Fix $g$ and $n$. We start by analyzing the complex
$S=\bigl(\FF_{\local^\Dual}\FF_\local\a\bigr)\(g,n\)$. The following lemma
identifies the underlying graded $\SS_n$-module $S^\sharp$.
\begin{lemma} \label{lemma.0}
a) There is an isomorphism
$$
S = \colim_{G\in\Iso\Gamma\(g,n\)} S(G) \cong \oplus_{G\in[\Gamma\(g,n\)]}
S(G)_{\Aut(G)} ,
$$
where
$$
S(G) = \bigoplus_{I\subset\Edge(G)} S(G,I) ,
$$
and
$$
S(G,I) = \local(G/I) \o \bigotimes_{v\in\VV(G/I)} \bigl(
\can(\pi_{G,I}^{-1}(v))^{-1} \o \local(\pi_{G,I}^{-1}(v)) \bigr) \o \a\(G\) .
$$

b) The hyperoperad structure of $\local$ induces a natural isomorphism
$S(G,I)\cong\Det(I)^{-1}\o\local(G)\o\a\(G\)$.
\end{lemma}
\begin{proof}
The proof of a) is similar to the formula of \ref{2.?} for $\MM^2\v$,
except that now, the cocycle factors associated to the two Feynman
transforms must be inserted at the appropriate points.

The proof of b) is as follows. The structure maps associated to the
morphism $\pi_{G,I}$ for the hyperoperads $\can$ and $\local$ induce
isomorphisms
\begin{gather*}
\local(G/I) \o \bigotimes_{v\in\VV(G/I)} \local(\pi_{G,I}^{-1}(v)) \cong
\local(G) , \\
\can(G/I) \o \bigotimes_{v\in\VV(G/I)} \can(\pi_{G,I}^{-1}(v)) \cong
\can(G) ,
\end{gather*}
and the ratio of these formulas gives
$$
\can(G/I)^{-1} \o \local(G/I) \o \bigotimes_{v\in\VV(G/I)} \bigl(
\can(\pi_{G,I}^{-1}(v))^{-1} \o \local(\pi_{G,I}^{-1}(v)) \bigr)
\cong \can(G)^{-1} \o \local(G) .
$$
Multiplying both sides by $\can(G/I)$ and observing that
$\can(G/I)\o\can(G)^{-1}\cong\Det(I)^{-1}$, part b) of the lemma follows.
\end{proof}

\newcommand{\Local}{\mathfrak{E}}

Recall that for any cocycle $\Local$ and any $\Local$-operad $\b$, the
differential in $\FF_\Local\b$ is a sum of two terms $\delta_{\b^*}+\p$,
where $\delta_{\b^*}$ is induced by the differential of $\b$, and $\p$ is
induced by the modular $\Local$-operad structure of $\b$, as explained in
\ref{FF}.

Applying this with $\Local=\local^\Dual$ and $\b=\FF_\local\a$, we find
that the differential in $\FF_{\local^\Dual}\FF_\local\a$ is a sum of three
terms $\delta_0+\delta_1+\delta_2$:
\begin{enumerate}
\item $\delta_0=\bigl(\p_0\bigr)_{\Aut(G)}:S(G)_{\Aut(G)}\to
S(G)_{\Aut(G)}$ is the map induced on $\Aut(G)$-coinvariants by
$\p_0:S(G,I)\to S(G,I)$, where $\p_0$ is the differential induced on the
summand $S(G,I)$ by the differential of $\a$;
\item the differential
$$
\delta_1 : S(G)_{\Aut(G)} \to \bigoplus_{\substack{H\in[\Gamma\(g,n\)]\\
|\Edge(H)|=|\Edge(G)|+1}} S(H)_{\Aut(H)}
$$
is induced by the differential $\p$ in $\FF_\local\a$, which itself is
induced by the modular $\local$-operad structure of $\a$;
\item $\delta_2=\bigl(\p_2\bigr)_{\Aut(G)}:S(G)_{\Aut(G)}\to
S(G)_{\Aut(G)}$ is the map induced on $\Aut(G)$-coinvariants by
$$
\p_2 : S(G,I) \to \bigoplus_{e\in I} S(G,I\setminus\{e\}) ,
$$
where $\p_2$ comes from the differential $\p$ on $\FF_{\local^\Dual}\b$,
and $\b=\FF_\local\a$.
\end{enumerate}
Thus, $\delta_0$ depends only on the internal differential of $\a$,
$\delta_1$ depends on the modular $\local$-operad structure of $\a$, while
$\delta_2$ is purely combinatorial and only depends on the graded stable
$\SS$-module structure underlying $\a$.

\begin{lemma} \label{lemma.1}
The map $\tau:(S,\delta_1+\delta_2)\to\a\(g,n\)$ is a weak equivalence of
complexes.
\end{lemma}

Let us first show how this lemma implies Theorem \ref{Feynman}. Observe
that $\delta_1$ has the effect of increasing both $|I|$ and $|\Edge(G)|$ by
$1$, while $\delta_2$ leaves $|\Edge(G)|$ unchanged, and decreases $|I|$ by
$1$. Therefore, $S$ is the total complex of a double complex
$(S_{\bull\bull},\delta_0,\delta_1+\delta_2)$, where
$$
S_{pq} = \bigoplus_{G\in[\Gamma\(g,n\)]} \biggl(
\bigoplus_{\substack{I\subset\Edge(G)\\|I|=q+2|\Edge(G)|}} S_{p+q}(G,I)
\biggr)_{\Aut(G)} ,
$$
where $S_{p+q}(G,I)$ is the degree $p+q$ subspace of the graded
$\SS_n$-module $S(G,I)$. Since there are a finite number of terms, indexed
by $G$ and $I$, contributing to this direct sum, this double complex has
$p+q\ge0$ and $q$ bounded below, and thus its associated spectral sequence
is convergent, yielding the desired implication.

We now turn to the proof of Lemma \ref{lemma.1}. Introduce the decreasing
filtration $\Phi$ of $S$ given by
$$
\Phi^q(S) = \bigoplus_{\substack{G\in[\Gamma\(g,n\)]\\|\Edge(G)|\ge q}}
S(G)_{\Aut(G)} .
$$
From the properties of $\delta_1$ and $\delta_2$ discussed above, we see
that
$$
\gr^\Phi(S) \cong (S,\delta_2) .
$$
This reduces the proof to that of the following lemma.
\begin{lemma} \label{lemma.2}
For each stable graph $G$ with $|\Edge(G)|>0$, the complex $(S(G),\p_2)$ is
acyclic.
\end{lemma}

Indeed, on taking $\Aut(G)$-coinvariants, this lemma implies that the
differential $\delta_2=(\p_2)_{\Aut(G)}$ on $S(G)_{\Aut(G)}$ is acyclic,
since the group $\Aut(G)$ is finite and we work over a field of
characteristic zero.

The proof of Lemma \ref{lemma.2} is based on the identification
\begin{equation} \label{equation.1}
(S(G),\p_2) \cong C_\bull(\Edge(G)) \o \local(G)^\sharp \o \a\(G\)^\sharp ,
\end{equation}
where, for any finite set $X$,
$$
C_\bull(X) = \bigoplus_{I\subset X} \Det(I)^{-1}
$$
is the augmented chain complex of the simplex with vertices $X$. Of course,
this chain complex is acyclic for $X$ non-empty. The identification
\ref{equation.1} follows from Lemma \ref{lemma.0} b).

This concluded the proof of Theorem \ref{Feynman}.
\end{proof}

\subsection{The Feynman transform and the cobar operad of a cyclic operad}
Let $\BB\a$ be the cobar operad the cyclic operad $\a$, introduced in
Section 3.2 of \cite{GK}. We may regard $\a$ as a modular operad
\ref{stab}; then $\Cyc(\FF\a)$ is related to $\BB\a$ by the formula
$$
\Cyc(\FF\a) \cong \Sigma\Susp\BB\a .
$$

\subsection{The Feynman transform and Vassiliev invariants}
\label{Vassiliev}
Vassiliev has introduced a filtered space $V=\bigcup_{k=0}^\infty V_k$ of
knot invariants of finite order (see Theorems 8 and 9 of \cite{BN}). The
associated graded space $W=\gr V$ is a commutative cocommutative Hopf
algebra. Let $P=\bigoplus P_k$ be its space of primitives. One of the chief
results of Kontsevich and Bar-Natan identifies $P_k$ with the lowest
homology groups of certain graph complex. In our language,
$$
P_k \cong \bigoplus_{\substack{k=g-1+n\\n>0}}
H_{1-g}(\Susp^{-1}\FF\Com)\(g,n\)_{\SS_n} ,
$$
where $\Com$ is the commutative operad.

\section{Modular operads and moduli spaces of curves}

In this section, we give some basic examples of modular operads, coming
from the theory of moduli spaces of stable algebraic curves. Throughout
this section, the base field is taken to be the field of complex numbers
$\C$.

\subsection{Orbifolds} Let $\CG$ be a groupoid in the category of varieties
over $\C$, with morphisms $\Mor(\CG)$, objects $\Ob(\CG)$, and source and
target maps $s,t:\Mor(\CG)\to\Ob(\CG)$. The groupoid $\CG$ is called
\begin{enumerate}
\item proper if the morphism $s\times
t:\Mor(\CG)\to\Ob(\CG)\times\Ob(\CG)$ is proper;
\item \'etale if $s$ and $t$ are \'etale;
\item smooth if $\Mor(\CG)$ and $\Ob(\CG)$ are smooth.
\end{enumerate}

An orbifold (smooth algebraic stack) is an equivalence class of smooth
proper \'etale group\-oids: two groupoids $\CG_1$ and $\CG_2$ are
equivalent if there is an \'etale map $f:\Ob(\CG_1)\to\Ob(\CG_2)$ and an
equivalence of categories $\CG_1\cong f^*\CG_2$, or more generally, if they
are joined by a chain of such equivalences. For more on orbifolds, see
\cite{DM} and \cite{Faltings}.

A sheaf $(\CS,p)$ on an orbifold $\CG$ is a sheaf $\CS$ on $\Ob(\CG)$
together with an isomorphism $p:s^*\CS\cong t^*\CS$. A global section of
such a sheaf is a global section $f$ of $\CS$ over $\Ob(\CG)$ such that
$s^*f$ and $t^*f$ are identifed by $p$.

The coarse space $|\CG|$ of an orbifold $\CG$ is the quotient of $\Ob(\CG)$
by the action of $\CG$, in other words, the space of isomorphism classes of
objects of $\CG$. Note that the coarse space $|\CG|$ need not be
smooth.

If $G$ is a group acting on an orbifold $\CG$, the quotient $\CG/G$ is the
orbifold with the same objects as $\CG$, and whose morphisms are
$G\times\Mor(\CG)$. The structure maps are defined as follows:
$$
s(g,x) = s(x) ,\quad t(g,x) = g(t(x)),\quad (g,x)\*(h,y)=(gh,h(x)\*y) .
$$
The coarse space of $\CG/G$ is isomorphic to the quotient $|\CG|/G$.

\subsection{Deligne-Mumford moduli spaces}\label{DM}
If $2(g-1)+n>0$, the (large) groupoid of smooth complex curves of genus $g$
with $n$ marked points, with isomorphisms as arrows, represents an orbifold
$\CM_{g,n}$, of dimension $3(g-1)+n$. As $g$ and $n$ are varied, we obtain
an $\SS$-orbifold, which we denote by $\CM$.

Knudsen \cite{Knudsen} proves that the (large) groupoid of stable complex
curves of genus $g$ with $n$ marked points, again with isomorphisms as
arrows, represents an orbifold $\Mbar_{g,n}$, of dimension $3(g-1)+n$.  As
$g$ and $n$ are varied, we obtain an $\SS$-orbifold, which we denote by
$\Mbar$, which contains $\CM$ as a dense open subset.

The dual graph $G(C,x_1,\dots,x_n)\in\Gamma\(g,n\)$ of a stable curve
$(C,x_1,\dots,x_n)\in\Mbar\(g,n\)$ is the labelled graph defined as
follows. Its flags are pairs $(K,y)$ where $y$ is either a nodal point or a
marked point $x_i$ and $K$ is a branch of the curve $C$ at $y$. (Note that
the curve has one branch at a marked point and two branches at a node.) Its
vertices are the components of $C$, its edges are the nodes, and its legs
are the points $x_i$. If $v\in G(C,x_1,\dots, x_n)$ is the vertex
corresponding to the component $K\in C$, label $v$ by the genus $g(v)$ of
the desingularization of $K$.

Given $G\in\Gamma\(g,n\)$, denote by $\CM_G\subset\Mbar\(g,n\)$ the
orbifold of stable curves whose dual graph is $G$; note that $\CM_G$ is
isomorphic to the orbifold $\CM\(G\)/\Aut(G)$. This gives a stratification
of $\Mbar\(g,n\)$ whose strata correspond to elements of $\Gamma\(g,n\)$;
the open stratum $\CM\(g,n\)$ corresponds to the graph with no edges. The
closure $\Mbar_G$ of $\CM_G$ is isomorphic to the orbifold
$\Mbar\(G\)/\Aut(G)$.

The $\SS$-orbifold $\Mbar$ is a modular operad $\Mbar$, with product
defined as follows: if $G\in\Gamma\(g,n\)$ is a stable graph, the
composition map
\begin{equation}\label{stratum}
\mu_G: \Mbar\(G\) = \prod_{v\in\VV(G)} \Mbar\(g(v),v\) \to\Mbar\(g,n\)
\end{equation}
is defined by gluing the marked points of the curves from
$\Mbar\(g(v),v\)$, $v\in\VV(g)$, according to the graph $G$ (see
\cite{GK}, 1.4.3). This map induces the embedding of $\Mbar_G$ as a
closed stratum of $\Mbar$.

Taking homology, we obtain a modular operad $H_\bull(\Mbar)$ in the
category of graded vector spaces. An algebra over this operad is the same
as a cohomological field theory in the sense of Kontsevich-Manin \cite{KM}.

\subsection{Differential forms with logarithmic singularities and
principal values}
Let $X$ be a compact $n$-dimensional complex manifold and $D\subset X$ a
divisor with normal crossings. Let $D^k\subset X$ be the locus of $k$-fold
self-intersection of $D$ (so that $D^0=X$ and $D^1=D$), with inclusion
morphisms $i^k:D^k\hookrightarrow X$ and $j^k:D^k\setminus
D^{k+1}\hookrightarrow X$. Let $\pi^k:\tilde D^k\to D^k\subset X$ be the
normalization morphism of $D^k$. The variety $\tilde{D}^k$ is smooth, and
the preimage
$$
\tildetilde{D}^{k+1} = (\pi^k)^{-1}(D^{k+1})
$$
is a divisor in $\tilde{D}^k$ with normal crossings.

Let $\CE^\bull_X$ be the complex of sheaves of $\Cinf$ differential forms
on $X$. The complex $\CE^\bull_X(\log D)$ of sheaves of $\Cinf$
differential forms with logarithmic singularities is the sheaf of
subalgebras of $j^1_*\CE^\bull_{X\setminus D}$ generated by $\CE^\bull_X$
and forms $df/f$ where $f$ is a holomorphic equation of $D$.

Let $\CE^\bull(X)$ and $\CE^\bull(X,\log D)$ be the spaces of global
sections of the sheaves $\CE^\bull_X$ and $\CE^\bull_X(\log D)$. Each of
the spaces $\CE^i(X,\log D)$ and $\CE^i(X)$ are nuclear Fr\'echet spaces,
since they are spaces of $\Cinf$ global sections of smooth vector bundles.

Let $\CC_{X,\bull}$ be the sheaf of de Rham currents on $X$. The space of
global sections $\CC_i(X)=\Gamma(X,\CC_{X,i})$ is the topological dual of
$\CE^i(X)$, and the differential $\delta$ on $\CC_{X,\bull}$ has degree
$-1$ and is adjoint to the exterior differential $d$ on $\CE^\bull_X$.

The principal value (Herrera-Lieberman \cite{HL}) is the continuous map of
graded sheaves
$$
\pv : \CE^\bull_X(\log D) \to \CC_{X,2n-\bull}
$$
defined as follows: if $U\subset X$ is an open set,
$\alpha\in\Gamma(U,\CE^i_X(\log D))$ and $\om\in\Gamma_c(U,\CE^{2n-i}_X)$,
$$
\< \pv\alpha , \om \> = \lim_{\eps\to0} \int_{|\phi|\ge\eps} \alpha \. \om ,
$$
where $\phi$ is a holomorphic defining equation of $D\cap U$ in $U$. (The
limit is independent of $\phi$.)

The Poincar\'e residue is the map of graded sheaves
$$
\Res : \CE^\bull_X(\log D) \to
\pi^1_*\CE^{\bull-1}_{\tilde{D}}(\tildetilde{D}^2)
$$
which measures the deviation of $\pv$ from being a map of complexes (Prop.\
5.3 of \cite{HL}):
\begin{equation} \label{residue}
\delta(\pv(\alpha)) - \pv(d\alpha) = \pv\bigl( 2\pi i\Res(\alpha) \bigr) ,
\quad \alpha\in\CE^\bull_X(\log D) .
\end{equation}

\subsection{Currents with logarithmic singularities}
Denote the image of the injective map
$$
\pi^k_*\CE^\bull_{\tilde{D}^k}(\log\tildetilde{D}^{k+1}) \xrightarrow{\pv}
\CC_{X,2(n-k)-\bull}
$$
by $\CC_{X,\bull}(D,k)$, and the sum of these spaces as $k$ varies between
$0$ and $k$ by $\CC_{X,\bull}(D)$. (Note that these spaces have zero
intersection, so the sum is direct.) We call elements of $\CC_{X,\bull}(D)$
currents with logarithmic singularities. By \ref{residue}, the differential
$\delta$ maps $\CC_{X,\bull}(D,k)$ to
$\CC_{X,\bull}(D,k)\oplus\CC_{X,\bull}(D,k+1)$; thus, $\CC_{X,\bull}(D)$ is
a complex of sheaves. Note that the spaces of global sections $\CC_i(X,D)$
are nuclear Fr\'echet spaces.
\begin{proposition} \label{poincare}
The inclusion $\CC_{X,\bull}(D)\hookrightarrow\CC_{X,\bull}$ is a weak
equivalence of complexes of sheaves.
\end{proposition}
\begin{proof}
We may assume that $X=\C^n$, with the divisor $D$ is given by the equation
$z_1\dots z_m=0$, $m\le n$; we must prove the weak equivalence for the
stalks of the two complexes of sheaves at $0\in\C^n$. If
$I\subset\{1,\dots,n\}$, let $\C^I\subset\C^n$ be the corresponding
coordinate subspace, let $D^I\subset\C^I$ be the divisor given by
$\prod_{i\in I}z_i=0$, and let $(\C^\times)^I=\C^I\setminus D^I$, with
embedding $j^I:(\C^\times)^I\hookrightarrow\C^n$. Let $i(I)$ be the least
element of the set $I$.

The graded sheaf $\CC_{\C^n,\bull}(D)$ decomposes as a direct sum
$$
\CC_{\C^n,\bull}(D) \cong
\bigoplus_{I\subset\{1,\dots,m\}} \CE^{2(n-|I|)-\bull}_{\C^I}(\log D^I) .
$$
We write the associated decomposition of a current $T$ lying in the stalk
$\CC_{\C^n,i}(D)_0$ of $\CC_{\C^n,i}(D)$ at $0$ as
$$
T = \sum_{I\subset\{1,\dots,m\}} j^I_*\pv(\alpha_I) ,
$$
where $\alpha_I$ is in the stalk at $0$ of $\CE^{2(n-|I|)-i}_{\C^I}(\log
D^I)$. Define a map $h:\CC_{\C^n,i}(D)_0\to\CC_{\C^n,i+1}(D)_0$ by
$$
h(T) = \frac{1}{2\pi i} \sum_{\substack{I\subset\{1,\dots,m\}\\|I|>0}}
j^{I\setminus\{i(I)\}}_*\pv\Bigl( \frac{dz_{i(I)}}{z_{i(I)}} \. \alpha_I
\Bigr) .
$$
It is easily seen that this is a contracting homotopy from
$\CC_{\C^n,\bull}(D)_0$ to $(\CC_{\C^n,\bull})_0$, proving the proposition.
\end{proof}

The above constructions generalize to the situation of a divisor $D$ with
normal crossings in an orbifold $X$. The orbifold $X$ may be represented by
a groupoid $\CG$, and the divisor $D$ gives rise to divisor with normal
crossings in $\Ob(\CG)$, invariant under the action of $\CG$ (that is,
$s^{-1}(D)=t^{-1}(D)$). The sheaves $\CE^\bull_\CG(\log D)$ and
$\CC_{\CG,\bull}(D)$ are defined to be the subsheaves of
$\CE^\bull_{\Ob(\CG)}(\log D)$ and $\CC_{\Ob(\CG),\bull}(D)$ invariant
under the action of $\CG$, that is, such that $s^*\om=t^*\om$. (Pullbacks
of currents on $\Ob(\CG)$ by the maps $s$ and $t$ are well-defined, since
these maps are \'etale.)

\subsection{The log-complex of $\Mbar_{g,n}$}
The compactification divisor $D_{g,n}=\Mbar_{g,n}\setminus\CM_{g,n}$ is a
divisor with normal crossings in the orbifold $\Mbar_{g,n}$, which
decomposes into an intersection of smooth divisors, corresponding to the
graphs in $\Gamma\(g,n\)$ with one edge. Let $\CC_\bull(\Mbar,D)$ and
$\CC_\bull(\Mbar)$ be the stable $\SS$-modules
$$
\CC_\bull(\Mbar,D)\(g,n\) = \CC_\bull(\Mbar_{g,n},D_{g,n})
\quad\text{and}\quad
\CC_\bull(\Mbar)\(g,n\) = \CC_\bull(\Mbar_{g,n}) .
$$
The operation of pushing forward currents along the inclusion of strata
makes these into modular operads of chain complexes, which are weakly
equivalent, and whose homology is the operad $H_\bull(\Mbar)$ of graded
vector spaces.

\subsection{The topological Feynman transform}
Recall \cite{Grothendieck} that nuclear Fr\'echet spaces form a symmetric
monoidal category $\nf$, with operation $\widehat{\o}$ (projective tensor
product). Furthermore, the opposite symmetric monoidal category $\nf^\op$
is identified, via the operation $V\mapsto V'$ (strong dual) with the
category $\df$ of nuclear DF-spaces, also with the projective tensor
product.

Let $C(\nf)$ and $C(\df)$ be the symmetric monoidal categories of bounded
chain complexes with finite-dimensional homology over $\nf$ and $\df$. The
strong dual identifies $C(\nf)^\op$ with $C(\df)$, and the homology of the
dual complex is naturally dual to the homology of the original complex.

Imitating the construction of the Feynman transform in the topological
setting, substituting the strong dual and projective tensor product for
their algebraic analogues, we obtain a functor $\FF^\TOP$, the topological
Feynman transform, from modular operads in $C(\nf)$ to modular
$\can$-operads in $C(\df)$. This functor has a homotopy inverse
$\FF_\om^\TOP$, constructed in the analogous way.

The stable $\SS$-module $\CC_\bull(\Mbar,D)$ is an example of a modular
operad in $C(\nf)$. By \ref{poincare}, its homology may be identified with
$H_\bull(\Mbar,D)$, the homology operad of the topological modular operad
$\Mbar$.

\subsection{The gravity operad}
Consider the stable $\SS$-module $\tGrav$, given by
$$
\tGrav\(g,n\) = \CE^\bull(\Mbar_{g,n},\log D_{g,n})' .
$$
For any graph $G\in\Gamma\(g,n\)$ with one edge, we have the adjoint of the
residue map
$$
\Res_G^* : \tGrav\(G\) \o \can(G)^{-1} \to \tGrav\(g,n\) .
$$
Iterating these maps, we may define $\Res_G^*$ for any stable graph.  The
maps $(2\pi i)^{|\Edge(G)|}\Res_G^*$ are the composition maps making
$\tGrav$ into a modular $\can^{-1}$-operad in $C(\df)$, which we call the
gravity operad.

Note that the homology $\Grav\(g,n\)$ of $\tGrav\(g,n\)$ form a modular
$\can^{-1}$-operad in the category of finite-dimensional graded vector
spaces, such that
$$
\Grav\(g,n\) \cong H_\bull(\CM_{g,n}) .
$$
The results of \cite{weil} show that the $S^1$-equivariant cohomology of a
topological conformal field theory in two dimensions is a modular algebra
over the suspension $\Susp\tGrav$. This paper also gives an explicit
presentation for the cyclic operad $\Cyc(\Sigma\Susp\tGrav)$.  (See also
\cite{gravity}.) This cyclic operad is formal, in the sense that there is a
weak equivalence between it and and its homology. (See \cite{GK} and
\cite{gravity}.) It seems unlikely that this is true for $\tGrav$ and its
homology $\Grav$.

Recall \ref{coboundary} the invertible stable $\SS$-module
$\Pp\(g,n\)=\Sigma^{-6(g-1)-2n}$, whose coboundary satisfies
$\local_\Pp\cong\can^2$. We see that $\Pp\tGrav$ is a modular $\can$-operad
in $C(\df)$.
\begin{proposition}
We have an isomorphism $\FF^\TOP_\can\Pp\tGrav\cong\CC_\bull(\Mbar,D)$.
\end{proposition}
\begin{proof}
We have identifications
\begin{align*}
\CC_\bull(\Mbar_{g,n},D_{g,n}) &\cong
\bigoplus_k \CC_\bull(\Mbar_{g,n},D_{g,n}) \\
&\cong \bigoplus_k \CE^{6(g-1)+2n-2k-\bull}(D_{g,n}^k,\log D_{g,n}^{k+1}) \\
&\cong \bigoplus_{G\in\Gamma\(g,n\)}
\CE^{6(g-1)+2n-2|\Edge(G)|-\bull}(\Mbar\(G\),D\(G\))^{\Aut(G)} ,
\end{align*}
where $D\(G\)=\Mbar\(G\)\setminus\CM\(G\)$. The component of the last sum
corresponding to $G$ is
$$
\MM\bigl( \Pp^{-1}\CE^\bull(\Mbar,D) \bigr) \(G\) ,
$$
and we obtain the sought after identification at the level of stable
$\SS$-modules. In fact, this identification also respects the compositions
of the two modular operads. It remains to check that the differentials
coincide; we leave this to the reader.
\end{proof}

As with any sort of cobar construction, the equation $\delta^2=0$ in the
Feynman transform $\FF^\TOP_\can\Pp\tGrav$ is precisely the associativity
of the composition in $\Pp\tGrav$. This gives a simple explanation of why
$\tGrav$ is a modular $\can^{-1}$-operad.

\section{Characteristics of cyclic operads} \label{Cyclic}

If $\v$ is a stable $\SS$-module, we can associate to it a symmetric
function $\Ch(\v)$, called its characteristic. In this section and the
next, we give formulas for $\Ch(\BB\a)$ in terms of $\Ch(\a)$, where $\a$
is a cyclic operad, and for $\Ch(\FF\a)$ in terms of $\Ch(\a)$, where $\a$
is a modular operad. The first of these formulas involves a generalization
of the Legendre transform, and the second a generalization of the Fourier
transform, from power series in one variable to symmetric functions in
infinitely many variables. Here, symmetric functions arise because of
well-known correspondence between the characters of the symmetric group and
the ring of symmetric functions. For further details on the theory of
symmetric functions, see Chapter 1 of Macdonald \cite{Macdonald}.

\subsection{Symmetric functions}

Consider the ring
$$
\Lambda = \varprojlim \Z\[x_1,\dots,x_k\]^{\SS_k}
$$
of symmetric functions (power series) in infinitely many variables. The
following standard symmetric functions
\begin{gather*}
h_n(x_i) = \sum_{i_1\le\dots\le i_n} x_{i_1}\dots x_{i_n} ,\quad
e_n(x_i) = \sum_{i_1<\dots<i_n} x_{i_1}\dots x_{i_n} , \\
p_n(x_i) = \sum_{i=1}^\infty x_i^n ,
\end{gather*}
are called respectively the complete symmetric functions, the elementary
symmetric functions and the power sums. It is a basic fact that
\begin{gather*}
\Lambda = \Z\[h_1,h_2,\dots\] = \Z\[e_1,e_2,\dots\] , \\
\Lambda\o\Q = \Q\[p_1,p_2,\dots\] ,
\end{gather*}
that is, that each of these three series of symmetric functions freely
generates $\Lambda$ (in the case of the power sums, over $\Q$). In
particular, $h_1=e_1=p_1$, while $h_2=\half(p_1^2+p_2)$ and
$e_2=\half(p_1^2-p_2)$.

Let $\sigma$ be an element of the symmetric group $\SS_n$, with cycles of
length $a_1\ge a_2\ge\dots\ge a_\ell$; thus $n=a_1+\dots+a_\ell$. The cycle
index of $\sigma$ is the symmetric function
$$
\psi(\sigma) = p_{a_1} \dots p_{a_\ell} \in \Lambda .
$$
The characteristic of a finite-dimensional $\SS_n$-module $V$ is the
symmetric function
$$
\ch_n(V) = \frac{1}{n!} \sum_{\sigma\in\SS_n} \Tr_V(\sigma)
\psi(\sigma) .
$$
It may be proved that $\ch_n(V)$ is in $\Lambda$, although it is only
evident from its definition that it is in $\Lambda\o\Q$.

We extend the definition of $\ch_n$ to graded $\SS_n$-modules by
$$
\ch_n(V) = \sum_i (-1)^i \ch_n(V_i) ,
$$
where $V_i$ is the degree $i$ component of $V$. Finally, the characteristic
of a graded $\SS$-module $\v=\{\v(n)\mid n\ge0\}$ such that $\v(n)$ is
finite-dimensional for all $n$ is
$$
\ch(\v) = \sum_{n=0}^\infty \ch_n(\v(n)) .
$$

We denote by $\rk:\Lambda\to\Q\[x\]$ the ring homomorphism which sends
$$
h_n \mapsto \frac{x^n}{n!} ,
$$
or equivalently, $p_1\mapsto x$ and $p_n\mapsto0$, $n>1$. If $V$ is an
$\SS_n$-module,
$$
\rk(\ch_n(V)) = \frac{\dim(V)x^n}{n!} .
$$
For this reason, we call $\rk$ the rank homomorphism.

\subsection{Plethysm} \label{plethysm}
Plethysm is the associative operation on $\Lambda$, denoted $f\circ g$,
characterized by the formulas
\begin{enumerate}
\item $(f_1+f_2)\circ g=f_1\circ g+f_2\circ g$;
\item $(f_1f_2)\circ g=(f_1\circ g)(f_2\circ g)$;
\item if $f=f(p_1,p_2,\dots)$, then $p_n\circ f=f(p_n,p_{2n},\dots)$.
\end{enumerate}
Note that under the rank homomorphism, plethysm is carried into composition
of power series.

There is a monoidal structure on the category of $\SS$-modules, with
tensor product
$$
(\v\circ\w)(n) = \bigoplus_{k=0}^\infty \Bigl( \v(k) \o
\bigoplus_{f:\{1,\dots,n\}\to\{1,\dots,k\}} \bigotimes_{i=1}^k
\w(f^{-1}(i)) \Bigr)_{\SS_k} .
$$
(An operad $\v$ is just an $\SS$-module with an associative composition
$\v\circ\v\to\v$.)

\begin{proposition} \label{composition}
$\ch(\v\circ\w) = \ch(\v)\circ\ch(\w)$
\end{proposition}

When $\v$ and $\w$ are ungraded, this is proved in Macdonald
\cite{Macdonald}. In the general case, the proof depends on an analysis of
the interplay between the minus signs in the Euler characteristic and the
action of symmetric groups on tensor powers of graded vector spaces.

\subsection{Characteristic of $\SS$-modules} \label{Ch:cyclic}
If $\v=\{\v\(n\)\mid n\ge1\}$ is a cyclic $\SS$-module, its characteristic
is
$$
\Ch(\v) = \sum_{n=1}^\infty \ch_n(\v\(n\)) .
$$
There is a forgetful functor from cyclic $\SS$-modules to $\SS$-modules,
obtained by restricting the action of $\v(n)=\v\(n+1\)$ from $\SS_{n\+}$ to
the subgroup $\SS_n$. The characteristics of $\v$ considered as a cyclic
$\SS$-module and an $\SS$-module are related by
\begin{equation} \label{Ch-ch}
\ch(\v) = \frac{\p\Ch(\v)}{\p p_1} .
\end{equation}

\subsection{Examples of characteristics}
To illustrate the above definitions, let us give some examples of
characteristics of cyclic operads.

\subsubsection{The commutative operad}\label{comm2}
For the commutative operad $\Com$, $\Com\(n\)$ is the trivial representation
of $\SS_n$ for all $n\ge3$, see \ref{comm}. It follows that
$\ch_n(\Com\(n\))=h_n$ for $n\ge3$, and hence that
$$
\Ch(\Com)
= \exp\Bigl(\sum_{n=1}^\infty \frac{p_n}{n} \Bigr) - \bigl( 1+h_1+h_2 \bigr) .
$$

\subsubsection{The associative operad}\label{ass2}
Since \ref{ass}
$
\Ass\(n\) \cong \Ind_{\Z_n}^{\SS_n} \k,
$
$$
\ch_n(\Ass\(n\)) = \sum_{d|n} \frac{\phi(d)}{n} p_d^{n/d} ,
$$
where $\phi(d)$ is the Euler function (the number of units in $\Z/d$.
Summing over $n\ge3$, we see that
$$
\Ch(\Ass)
= - \sum_{n=1}^\infty \frac{\phi(n)}{n} \log(1-p_n) - (h_1 + h_2) .
$$

\subsubsection{The Lie operad}\label{lie2}
In \ref{Lie} below, we will prove that the characteristic of the 
Lie operad is
$$
\Ch(\Lie)
= (1-p_1) \sum_{n=1}^\infty \frac{\mu(n)}{n} \log(1-p_n) + p_1 ,
$$
where $\mu(n)$ is the M\"obius function.

\subsection{The Legendre transform}
Classically, the Legendre transform of a convex function $f:\R\to\R$ is
the function
$$
(\CL f)(\xi) = g(\xi) = \max_x \bigl( x\xi-f(x) \bigr) .
$$
(See Section 3.3 of Arnold \cite{Arnold}.) Setting $\xi=f'(x)$,
we see that
\begin{equation} \label{legendre}
g \circ f' + f = xf' .
\end{equation}
Suppose that, instead of being a convex function, $f(x)$ is a formal power
series of the form
\begin{equation} \label{x-squared}
f(x) = \sum_{n=2}^\infty \frac{a_nx^n}{n!} \in \Q\[x\] ,
\end{equation}
where $a_2\ne0$; we denote this set of power series by $\Q\[x\]_*$.
The equation \ref{legendre} defines a unique power series
$(\CL f)(\xi)=g(\xi)\in\Q\[\xi\]_*$, which we again call the Legendre
transform.
\begin{proposition} \label{inverse}
If $f$ and $g$ are series of the form \ref{x-squared}, then $\CL f=g$ if
and only if $f'$ and $g'$ are inverse under composition, that is,
$$
g' \circ f' = x .
$$
\end{proposition}
\begin{proof}
Taking the derivative of \ref{legendre}, we see that
$$
(g'\circ f')f'' + f' = xf'' + f' .
$$
Cancelling $f'$ from each side and dividing by $f''$, which is invertible
in $\Q\[x\]$ by hypothesis, we find that $g'\circ f'=x$. The same reasoning
proves the converse.
\end{proof}

As  a consequence of this proposition, we see that $\CL$ is involutive:
$\CL(\CL f)=f$.

\subsection{$\CL$ and trees}
Let $\v$ be a cyclic $\SS$-module with $\v\(n\)=0$ for $n\le2$. The cyclic
$\SS$-module $\TT\v$ was defined in \ref{TT}.
\begin{proposition}
Let $a_n=\chi(\v\(n\))$ and $b_n=\chi(\TT\v\(n\))$ be the Euler
characteristics of the components of $\v$ and $\TT\v$. If
$$
f(x) = \frac{x^2}{2} - \sum_{n=3}^\infty \frac{a_nx^n}{n!} \text{ and }
g(x) = \frac{x^2}{2} + \sum_{n=3}^\infty \frac{b_nx^n}{n!} ,
$$
then $g=\CL f$.
\end{proposition}
\begin{proof}
It is a corollary of Theorem 3.3.2 of \cite{GK} that $g'\circ f'=x$. The
results follows by \ref{inverse}.
\end{proof}

With the notation of the proposition,
\begin{equation} \label{sum-over-trees}
b_n = \sum_{\text{$n$-trees $T$}} \prod_{v\in\VV(T)} a_{n(v)} .
\end{equation}
In fact, the proposition remains true for an arbitrary sequence of rational
numbers $\{a_3,a_4,\dots\}$, if we define $\{b_3,b_4,\dots\}$ by
\ref{sum-over-trees}.

\subsection{The Legendre transform for symmetric functions}
Denote by $\Lambda_*$ the set of symmetric functions such that
$\rk(f)\in\Q\[x\]_*$.
\begin{theorem}
(a) If $f\in\Lambda_*$, there is a unique element $g=\CL f\in\Lambda_*$
such that
\begin{equation} \label{Legendre}
g\circ\frac{\p f}{\p p_1} + f = p_1 \frac{\p f}{\p p_1} .
\end{equation}
We call $\CL:\Lambda_*\to\Lambda_*$ the Legendre transform.

\noindent(b) The Legendre transform of symmetric functions is compatible
with that of power series, in the sense that the following diagram
commutes:
$$\begin{diagram}
\node{\Lambda_*} \arrow{e,t}{\CL} \arrow{s,l}{\rk}\node{\Lambda_*}
\arrow{s,l}{\rk} \\
\node{\Q\[x\]_*} \arrow{e,t}{\CL} \node{\Q\[x\]_*}
\end{diagram}$$

\noindent(c) The symmetric functions
$$
\frac{\p(\CL f)}{\p p_1} \text{ and } \frac{\p f}{\p p_1}
$$
are plethystic inverses. (Note that, unlike for power series, this equation
does not determine $\CL f$.)

\noindent(d) The transformation $\CL$ is an involution, that is,
$\CL\CL=\Id$.
\end{theorem}
\begin{proof}
If $f\in\Lambda_*$, then $\p f/\p p_1$ is invertible with respect to
plethysm. Thus \ref{Legendre} defines $g\in\Lambda_*$ uniquely, proving
(a). Part (b) is obvious, since $\rk$ transforms plethysm into composition.

In proving (c), we need an analogue of the chain rule for $\p/\p p_1$ acting
on $\Lambda$:
$$
\frac{\p}{\p p_1} (u\circ v) = \Bigl( \frac{\p u}{\p p_1} \circ v \Bigr)
\frac{\p v}{\p p_1} .
$$
This formula is proved by checking that both sides are compatible with the
rules (1-3) defining plethysm \ref{plethysm}.

Using this, the reasoning needed to prove (c) is formally identical to that
in the proof of \ref{inverse}.

To prove (d), we note that (c) implies
$$
p_1 \frac{\p f}{\p p_1} = \Bigl( p_1 \frac{\p g}{\p p_1} \Bigr) \circ
\frac{\p f}{\p p_1} .
$$
This shows that
$$
g\circ\frac{\p f}{\p p_1}\circ\frac{\p g}{\p p_1} + f\circ\frac{\p g}{\p p_1}
= \Bigl( p_1 \frac{\p f}{\p p_1} \Bigr) \circ \frac{\p g}{\p p_1}
= \Bigl( p_1 \frac{\p g}{\p p_1} \Bigr) \circ \frac{\p f}{\p p_1}
\circ\frac{\p g}{\p p_1} .
$$
Cancellation proves that
$$
g + f\circ\frac{\p g}{\p p_1} = p_1 \frac{\p g}{\p p_1} ,
$$
and hence that $f=\CL g$.
\end{proof}

For example, $\CL h_2=e_2$ and vice versa.

The following theorem is related to results of Otter \cite{Otter} and
Hanlon-Robinson \cite{HR} on the enumeration of unrooted trees.
\begin{theorem} \label{7.17}
Let $\v$ be a cyclic $\SS$-module such that $\v\(n\)=0$ for $n\le2$ and
$\v\(n\)$ is finite dimensional for all $n$. Define the elements of
$\Lambda_*$
$$
\text{$f = e_2 - \Ch(\v)$ and $g = h_2 + \Ch(\TT\v)$.}
$$
Then $g=\CL f$.
\end{theorem}
\begin{proof}
By definition of $\CL$, we must prove that
$$
\bigl( h_2 + \Ch(\TT\v) \bigr) \circ \bigl( p_1 - \ch(\v) \bigr)
+ \bigl( e_2 - \Ch(\v) \bigr) = p_1 \bigl( p_1 - \ch(\v) \bigr) ,
$$
since by \ref{Ch-ch}, $f'=p_1-\ch(\v)$. A little rearrangement shows that
this formula is equivalent to
\begin{equation} \label{Otter}
\Ch(\TT\v) \circ \bigl( p_1 - \ch(\v) \bigr)
= \Ch(\v) - e_2\circ\ch(\v) .
\end{equation}
Indeed, the formulas $h_2\circ(a+b)=h_2\circ a+h_2\circ b+ab$ and
$h_2\circ(-a)=e_2\circ a$ show that
$$
h_2 \circ \bigl( p_1 - \ch(\v) \bigr) = h_2 + e_2\circ\ch(\v) - p_1\ch(\v) .
$$
Equation \ref{Otter} now follows from the formula $p_1^2=h_2+e_2$.

We prove \ref{Otter} by constructing a differential graded $\SS$-module
$\c=\{\c(n)\}$ such that the left-hand side of \ref{Otter} equals
$\ch(\c)$, and the right-hand side equals $\ch(H_\bull(\c))$. Define the
$\SS$-module underlying $\c$ to be the plethysm $\x\circ\w$, where the
$\SS$-modules $\x$ and $\w$ are defined by
\begin{align*}
\x(n) &= \begin{cases} 0 , & n\le2 , \\ (\TT\v)\(n\) , & n\ge3 ,
\end{cases} \\
\w(n) &= \begin{cases} 0 , & n=0 , \\ \k , & n=1 , \\
\Sigma\Res^{\SS_{n\+}}_{\SS_n}\v\(n+1\) , & n\ge2 . \end{cases}
\end{align*}
(Here, $\Sigma$ is the suspension functor on graded $\SS_n$-modules.) It
follows from \ref{composition} that $\ch(\c)$ equals the right-hand side of
\ref{Otter}.

We now construct a differential $\delta$ on $(\x\circ\w)(n)$. We say that a
vertex $v$ of a tree $T$ is a boundary vertex if exactly one of its flags
forms part of an edge; denote by $\beta(T)\subset\VV(T)$ the set of
boundary vertices of $T$. Then
\begin{equation} \label{coloured}
(\x\circ\w)(n) = \bigoplus_{\text{$n$-trees T}} \bigoplus_{B\subset\beta(T)}
\Bigl( \bigotimes_{v\in\VV(T)\setminus B} \v\(v\) \o
\bigotimes_{v\in B} \Sigma\v\(v\) \Bigr) .
\end{equation}
On the summand of \ref{coloured} associated to $(T,B)$, define the
differential
$$
\delta = \sum_{v\in B} \delta_v ,
$$
where $\delta_v$ is the natural identification, of degree $-1$, between
this summand and the summand associated to $(T,B\setminus\{v\})$.

Clearly $(\x\circ\w)(n)$ splits into a sum of subcomplexes $\c_T$ indexed
by the $n$-trees $T$. If $T$ has at least one non-boundary vertex, the
complex $\c_T$ is isomorphic to the tensor product of the graded vector
space $\v\(T\)$ and the augmented chain complex of the simplex with
vertices $\beta(T)$, and is thus contractible. As observed by Jordan
\cite{Jordan}, there remain trees with either one vertex or one edge. We
consider each of these cases separately.
\begin{enumerate}
\item The characteristic $\ch(\c_T)$ summed over trees $T$ with one vertex
equals $\Ch(\v)$.
\item The characteristic $\ch(\c_T)$ summed over trees $T$ with one edge
has two contributions: the terms with $B$ empty, which sum to
$h_2\circ\ch(\v)$, and the terms in which $|B|=1$, which sum to
$-\ch(\v)^2$. The sum of these two terms is $-e_2\circ\ch(\v)$.
\qed\end{enumerate}
\def\qed{}
\end{proof}

\subsection{The involution $\tilde\omega$ and the characteristic of the cobar
operad} Using the theorem just proved, we now write a formula for
$\Ch(\BB\a)$, where $\a$ is a cyclic operad. Up to differential, $\BB\a$ is
the cyclic operad $\TT\Susp^{-1}\Sigma^{-1}\a^*$, and thus
$\Ch(\BB\a)=\Ch(\TT\Susp^{-1}\Sigma^{-1}\a^*)$. Since
$\Ch(\Sigma^{-1}\a^*)=-\Ch(\a)$, it suffices to determine the effect of
$\Susp$ and $\Sigma$ on $\Ch(\v)$.

Denote by $\om:\Lambda\to\Lambda$ the ring homomorphism such that
$\om(h_n)=e_n$, $n\ge1$. If $V$ is a finite-dimensional $\SS_n$-module,
$$
\om(\ch_n(V)) = \ch_n(\sgn_n\o V) ,
$$
and thus $\om$ is an involution. Note also that $\om(p_n)=(-1)^{n-1}p_n$.

We also need a modified involution $\tom$, defined by
$\tom(h_n)=(-1)^ne_n$, or equivalently $\tom(p_n)=-p_n$. Thus, if $\v$ is a
cyclic $\SS$-module such that $\v\(n\)$ is finite-dimensional for each $n$,
\begin{equation}
\Ch(\Susp\v) = \tom(\Ch(\v)) .
\end{equation}
\begin{corollary} \label{BB}
Let $\a$ be a cyclic operad such that $\a\(n\)=0$ for $n\le2$ and $\a\(n\)$
is finite-dimensional for each $n$, and let $\BB\a$ be its cobar operad.
Then
$$
h_2 + \Ch(\BB\a) = \CL\tom(h_2+\Ch(\a)) .
$$
\end{corollary}

Recall \cite{GK} that $\BB\BB\a$ is weakly equivalent to $\a$, which
suggests that the transform $\CL\tom:\Lambda_*\to\Lambda_*$ should be an
involution. This follows from the next result.
\begin{proposition}
If $f\in\Lambda_*$, $-\CL\tom f=\CL(-f)$.
\end{proposition}
\begin{proof}
By Ex.\ 8.1 of Macdonald \cite{Macdonald}, if $u$ and $v$ are symmetric
functions, $u\circ(-v)=(\tom u)\circ v$. If $g=\CL(\tom f)$, we see
that the defining equation
$$
(\tom f)\circ \frac{\p g}{\p p_1} + g = p_1 \frac{\p g}{\p p_1}
$$
is equivalent to
$$
- f\circ \frac{\p(-g)}{\p p_1} + (-g) = - p_1 \frac{\p(-g)}{\p p_1} ,
$$
that is, $-g=\CL(-f)$.
\end{proof}

\subsection{Example: the Lie operad} \label{Lie}
The Lie operad $\Lie$ is weakly equivalent to the cobar operad $\BB\Com$ of
the commutative operad,  see \cite{GK}, and thus $h_2+\Ch(\Lie)$ is 
the Legendre transform
of
\begin{align*}
\tom(h_2+\Ch(\Lie)) &= \tom\Bigl(
\exp\Bigl(\sum_{n=1}^\infty \frac{p_n}{n} \Bigr) - \bigl( 1+h_1 \bigr)
\Bigr) \\
&= \exp\Bigl(\sum_{n=1}^\infty - \frac{p_n}{n} \Bigr) - \bigl( 1-p_1 \bigr) .
\end{align*}
The Legendre transform of this symmetric function is
$$
(1-p_1) \sum_{n=1}^\infty \frac{\mu(n)}{n} \log(1-p_n) + p_1 .
$$
It follows that
$$
\Ch(\Lie) = (1-p_1) \sum_{n=1}^\infty \frac{\mu(n)}{n} \log(1-p_n) + h_1 -
h_2,
$$
as was promised in \ref{lie2}.

\section{Characteristics of modular operads} \label{Modular}

\subsection{The ring $\Lambda\(\hbar\)$}\label{lhbar}
Consider the ring $\Lambda\(\hbar\)$ of Laurent series with coefficients in
$\Lambda$. This ring has a descending filtration
$$
F^m\Lambda\(\hbar\) = \Bigl\{ \sum f_i\hbar^i \,\big|\, f_i \in
F^{m-2i}\Lambda \Bigr\} ,
$$
inducing a topology on $\Lambda\(\hbar\)$. If $f\in\Lambda$, the plethysm
$f\circ(-):\Lambda\to\Lambda$ extends to $\Lambda\(\hbar\)$ by retaining
axioms (1) and (2) of \ref{plethysm} and replacing (3) by

($3'$) $p_n \circ f(\hbar,p_1,p_2,\dots) =
f(\hbar^n,p_n,p_{2n},\dots)$.

\subsection{The characteristic of a stable $\SS$-module}
The characteristic of a stable $\SS$-module $\v$ is the element of
$\Lambda\(\hbar\)$ given by the formula
$$
\CCh(\v) = \sum_{2(g-1)+n>0} \hbar^{g-1} \ch_n(\v\(g,n\)) .
$$
The stability condition ensures that $\CCh(\v)\in F^1\Lambda\(\hbar\)$.
If $\v$ is a cyclic $\SS$-module, the characteristic $\CCh(\v)$ of $\v$
considered as a stable $\SS$-module, equals $\hbar^{-1}\Ch(\v)$.

Note that as in the case of cyclic $\SS$-modules, we have
\begin{equation} \label{tom}
\CCh(\Susp\v) = \tom(\CCh(\v)) .
\end{equation}

Our goal in this section is to present formulas for $\CCh(\MM\a)$ and
$\CCh(\MM_{\Det(\Edge)}\a)$ in terms of $\CCh(\v)$. This will also
permit us to give formulas for $\CCh(\FF\a)$ and $\CCh(\FF^{-1}\a)$.

\subsection{Plethystic exponential}
For $f\in F^1\Lambda\(\hbar\)$, let
$$
\Exp(f) = \Bigl( \sum_{n=0}^\infty h_n \Bigr) \circ f
= \exp\Bigl(\sum_{n=1}^\infty \frac{p_n}{n} \Bigr) \circ f .
$$
Note that
$$
\Exp(f+g) = \Exp(f)\Exp(g) ,
$$
and that under specialization $\rk:\Lambda\(\hbar\)\to\Q\[x\]\(\hbar\)$,
the map $\Exp$ goes into exponentiation
$$
f(\hbar,x) \mapsto e^{f(\hbar,x)} .
$$

\begin{proposition}
If $\v$ is a stable $\SS$-module, let $\Exp_n(\v)$ be the stable
$\SS$-module such that
$$
\Exp_n(\v)\(g,n\) = \biggl(
\bigoplus_{\substack{f:I\to\{1,\dots,n\} \\ g_1+\dots+g_n=g}}
\Ind_{\Aut(f)}^{\Aut(I)}
\Bigl( \bigotimes_{i=1}^n \v\(g_i,f^{-1}(i)\) \Bigr) \biggr)_{\SS_n} ,
$$
where $\Aut(f)=\Aut(f^{-1}(1)) \times \dots \times \Aut(f^{-1}(n))$. Then
$$
\Exp(\CCh(\v)) = \sum_{n=0}^\infty \hbar^{-n} \CCh(\Exp_n(\v)) .
$$
\end{proposition} 
\begin{proof}
This follows from \ref{composition} and the definition of $\Exp(f)$,
$f\in\Lambda\(\hbar\)$.
\end{proof}

Informally, the stable $\SS$-module $\Exp_n(\v)$ may be thought of as
representing disconnected graphs with $n$ vertices and no edges: all of
its flags are legs.

The following proposition is essentially due to Cadogan \cite{Cadogan},
although he does not use the notation $\Exp$.
\begin{proposition}
The map $\Exp:F^1\Lambda\(\hbar\)\to1+F^1\Lambda\(\hbar\)$ is invertible
over $\Q$, with inverse
$$
\Log(f) = \sum_{n=1}^\infty \frac{\mu(n)}{n} \log(p_n) \circ f
= \sum_{n=1}^\infty \frac{\mu(n)}{n} \log(p_n\circ f) .
$$
\end{proposition}
\begin{proof}
\begin{align*}
\Log(\Exp(f)) &= \sum_{n=1}^\infty \frac{\mu(n)}{n} \log(p_n)
\circ \exp\Bigl( \sum_{m=1}^\infty \frac{p_m}{m} \Bigr) \circ f \\
&= \sum_{n=1}^\infty \frac{\mu(n)}{n} \log
\exp\Bigl( \sum_{m=1}^\infty \frac{p_{nm}}{m} \Bigr) \circ f \\
&= \sum_{n=1}^\infty \sum_{d|n} \frac{\mu(d)p_n\circ f}{n} = f .
\qed\end{align*}
\def\qed{}
\end{proof}

\subsection{The inner product on $\Lambda$}
To a partition $\lambda=(1^{m_1}2^{m_2}\dots)$, where $m_k=0$ for $k\gg0$,
is associated a monomial
$$
p_\lambda = p_1^{m_1} p_2^{m_2} \dots .
$$
These monomials form a topological basis of $\Lambda$. Let $\Lambda_\alg$ be
the space of finite linear combinations of the $p_\lambda$. The standard
inner product on $\Lambda_\alg$ is determined by the formula
$$
\< p_\lambda , p_\mu \> = \delta_{\lambda\mu} \prod_{i=1}^\infty i^{m_i} m_i!
$$
Note in particular that $\<p_i,p_j\>=i\delta_{ij}$; the inner product on
$\Lambda_\alg$ is the standard extension of the inner product on a vector
space to its symmetric algebra (Fock space).
\begin{proposition}
If $V$ and $W$ are $\SS_n$-modules,
$$
\< \ch_n(V) , \ch_n(W) \> = \dim \Hom_{\SS_n}(V,W) .
$$
\end{proposition}
\begin{proof}
This statement is well-known in the theory of symmetric functions: it
follows from the fact that the Schur functions form an orthonormal basis of
$\Lambda_\alg$. 
\end{proof}

We extend the inner product on $\Lambda_\alg$ to a $\Q\(\hbar\)$-valued
inner product on $\Lambda_\alg\(\hbar\)$ by $\Q\(\hbar\)$-bilinearity. If
$f\in\Lambda_\alg\(\hbar\)$, let
$D(f):\Lambda\(\hbar\)\to\Lambda\(\hbar\)$ be the adjoint of
multiplication by $f$ with respect to this inner product. The following
proposition is Ex.\ 5.3 of Macdonald \cite{Macdonald}.
\begin{proposition} \label{D(f)}
If $f=f(\hbar,p_1,p_2,\dots)\in\Lambda_\alg\(\hbar\)$, then
$$\textstyle
D(f) = f\bigl( \hbar , \frac{\p}{\p p_1} , 2\frac{\p}{\p p_2} ,
3\frac{\p}{\p p_3} , \dots \bigr) .
$$
\end{proposition}

\begin{proposition}
Let $k\le n$, $V$ be an $\SS_k$-module, and $W$ be an $\SS_n$-module. Then
$$
D(\ch_k(V)) \ch_n(W) = \ch_{n-k} \Hom_{\SS_k} \Bigl( V ,
\Res^{\SS_n}_{\SS_k\times\SS_{n-k}} W \Bigr) .
$$
\end{proposition}
\begin{proof}
This follows by taking adjoints on both sides of the formula
$$
\ch_j(U) \ch_k(V) = \ch_{j+k} \Ind^{\SS_{j+k}}_{\SS_j\times\SS_k} (U\o V) .
\qed$$
\def\qed{}
\end{proof}

\subsection{A Laplacian on $\Lambda\(\hbar\)$}
We now introduce an analogue of the Laplacian on $\Lambda\(\hbar\)$, given
by the formula
$$
\Delta = \sum_{n=1}^\infty \hbar^n
\left( \frac{n}{2} \frac{\p^2}{\p p_n^2} + \frac{\p}{\p p_{2n}} \right) .
$$
Note that $\Delta$ is homogeneous of degree zero, and thus preserves the
filtration of $\Lambda\(\hbar\)$. Under specialization
$\rk:\Lambda\(\hbar\)\to\Q\[x\]\(\hbar\)$, the operator $\Delta$
corresponds to the Laplacian $\frac{\hbar}{2}\frac{d}{dx^2}$ on the line.
\begin{proposition}
$D(\Exp(\hbar h_2))=\exp(\Delta)$
\end{proposition}
\begin{proof}
By \ref{D(f)}, it suffices to substitute $n\p/\p p_n$ for $p_n$ on the
right-hand side of
$$
\Exp(\hbar h_2) = \exp\Bigl(\sum_{n=1}^\infty \frac{p_n}{n} \Bigr)
\circ \Bigl( \frac{\hbar}{2} (p_1^2+p_2) \Bigr)
= \exp\Bigl(\sum_{n=1}^\infty \frac{\hbar^n}{2n} (p_n^2+p_{2n}) \Bigr) .
$$
\end{proof}

\begin{theorem}\label{modular-character}
If $\v$ is a stable $\SS$-module, then
$$
\CCh(\MM\v) = \Log\bigl( \exp(\Delta)\Exp(\CCh(\v)) \bigr) .
$$
\end{theorem}
\begin{proof}
We start by neglecting $\hbar$ and explain the appearance of the sum over
graphs on the right-hand side of the formula. Formally, we set $\hbar=1$;
this is legitimate if $\v\(g,n\)=0$ for $g\gg0$.

Applying $\Exp$ to $\CCh(\v)$, we obtain the stable $\SS$-module
representing possibly disconnected graphs each component of which has one
vertex. Applying $D(h_2)$ to $\Exp(\CCh(\v))$ gives the sum over all ways of
joining two legs (or flags) of such a graph; $h_2$ arises because the two
ends of an edge are indistinguishable. (If edges carried a direction,
we would replace $h_2$ by $p_1^2$, the characteristic of the regular
representation of $\SS_2$.)

Similarly, applying $D(\Exp(\hbar h_2))$ to $\Exp(\CCh(\v))$ gives the sum
over all ways of joining together any number $N$ of pairs of legs by
edges. In this way, we see (recall that $\hbar$ temporarily equals
$1$) that
$$
\exp(\Delta) \Exp(\CCh(\v)) = \CCh(\w) ,
$$
where $\w$ is the stable $\SS$-module such that
\begin{equation} \label{W}
\w\(g,n\) = \bigoplus_G \v\(G\)_{\Aut(G)} ,
\end{equation}
where $G$ runs over all possibly disconnected, labelled $n$-graphs such
that each component is stable. But $\w=\Exp(\MM\v)$, since $\MM\v$ is
defined in a similar way, but summing only over connected graphs.

To finish the proof, we must account for the powers of $\hbar$ in each term
of \ref{W}. Each term $\ch_n(\v\(g,n\))$ in $\CCh(\v)$ comes with a factor
of $\hbar^{g-1}$. The term of $\Exp(\CCh(\v))$ corresponding to a labelled
graph $G$ (with each component having one vertex) comes with a factor of
$\hbar$ raised to the power
$$
- |\VV(G)| + \sum_{v\in\VV(G)} g(v) .
$$
Each new edge introduced by the action of $D(\Exp(\hbar h_2))$ contributes
a factor of $\hbar$. Therefore, the term in \ref{W} corresponding to a
labelled graph $G$ comes with a factor of $\hbar$ raised to the power
$$
- \chi(G) + \sum_{v\in\VV(G)} g(v) .
$$
Applying $\Log$ has the effect of discarding all the disconnected graphs
$G$. If $G$ is connected, the power of $\hbar$ in question equals $g(G)-1$,
where $g(G)$ is defined in \ref{g(G)}.
\end{proof}

Recall that $\FF$ denotes the Feynman transform $\FF_\1$ associated to the
trivial cocycle.
\begin{corollary} \label{ff}
$\CCh(\FF\a) = \Log\bigl( \exp(-\Delta))\Exp(\CCh(\a)) \bigr)$
\end{corollary}
\begin{proof}
If $\a$ is a modular $\can$-operad, the stable $\SS$-modules $\MM\a$ and
$\FF^{-1}\a$ have the same characteristic, showing that
$$
\CCh(\FF^{-1}\a) = \Log\bigl( \exp(\Delta))\Exp(\CCh(\a)) \bigr) .
$$
Since $\FF$ is a homotopy inverse of $\FF^{-1}$, and characteristics are
homotopy invariant, the result follows.
\end{proof}

Note that \ref{ff} may also be proved by calculating $\CCh(\MM_\can\v)$;
this is given by the stated formula, since the effect of twisting by
$\Det(\Edge)$ is to attach a suspension to each edge, which changes the
operator $D(\Exp(\hbar h_2))$ into $D(\Exp(-\hbar h_2))=\Exp(-\Delta)$.

\subsection{Plethystic Fourier transform}
Let us give a formal interpretation of the previous theorem in terms of the
Fourier transform on the infinite-dimensional vector space
$\Spec(\Lambda_\R)\cong\R^\infty$, with coordinates $p_1,p_2,\dots$, where
$\Lambda_\R=\Lambda_\alg\o\R$. This space has a translation invariant
Riemannian metric
\begin{equation} \label{metric}
\<p_i,p_j\> = i\delta_{ij} .
\end{equation}
We denote the function $p_n\o1\in\Lambda_\R\o\Lambda_\R$ by $p_n$, and the
function $1\o p_n$ by $q_n$.

Let $d\mu$ be the formal Gaussian measure
$$
d\mu = \prod_{\text{$n$ odd}} e^{-p_n^2/2n\hbar^n}
\prod_{\text{$n$ even}} e^{-p_n^2/2n\hbar^n+p_n/n\hbar^{n/2}}
\frac{dp_n}{e^{c(n)/2n}\sqrt{2\pi n \hbar^n}}
$$
on $\Spec(\Lambda_\R)$, where $c(n)=\half((-1)^n+1)$. We may rewrite this
measure as
$$
d\mu = \Exp(-e_2/\hbar) \, \prod_{n=1}^\infty
\frac{dp_n}{e^{c(n)/2n}\sqrt{2\pi n\hbar^n}} ;
$$
it is the translate of the Gaussian measure associated to the metric
\ref{metric} by the vector
$(p_1,p_2,\dots)=(0,\hbar,0,\hbar^2,0,\hbar^3,0,\dots)$.

\newcommand{\hint}{\int^*}

If $p^\alpha$ is a monomial in the variables $p=(p_1,p_2,\dots)$, where
$\alpha$ is a multi-index $(\alpha_1,\alpha_2,\dots)$, define a
power series
$$
\hint_{\R^\infty} p^\alpha \, d\mu(p) \in \Z\[\hbar\]
$$
by the formula
\begin{align*}
\hint_{\R^\infty} p^\alpha \, d\mu(p) &=
\prod_{\text{$n$ odd}} \int_{-\infty}^\infty p_n^{\alpha_n}
e^{-p_n^2/2n\hbar^n} \frac{dp_n}{\sqrt{2\pi n\hbar^n}} \\
& \quad \times \prod_{\text{$n$ even}} \int_{-\infty}^\infty p_n^{\alpha_n}
e^{-p_n^2/2n\hbar^n+p_n/n\hbar^{n/2}} \frac{dp_n}{e^{1/2n}\sqrt{2\pi
n\hbar^n}} .
\end{align*}
This formula makes sense because almost all terms equal $1$, and it may be
defined in purely algebraic way by induction on $|\alpha|$: for $\alpha=0$
it equals $1$, and the induction step is performed by means of integration
by parts in one of the variables $p_n$. Extend the operation
$f\mapsto\hint_{\R^\infty}f\,d\mu(p)$ to a map from $\Lambda\(\hbar\)$ to
$\Z\(\hbar\)$ by linearity.

We may now restate \ref{modular-character} in the form of a (Gaussian)
Fourier transform. In the course of the proof, we make use of another
formal integral $\hint_{\R^\infty}f\,d\nu(p)$, whose definition is similar
to $\hint_{\R^\infty}f\,d\mu(p)$ and is given by the formula
$$
\hint_{\R^\infty} p^\alpha \, d\nu(p) = \prod_n \int_{-\infty}^\infty
p_n^{\alpha_n} e^{-p_n^2/2n\hbar^n} \frac{dp_n}{\sqrt{2\pi n\hbar^n}}
$$
\begin{theorem} \label{Wick}
As a function of $\hbar$ and $q=(q_1,q_2,\dots)$,
$$
\hbar^{-1}h_2 + \CCh(\MM\v) = \Log \hint_{\R^\infty}
\Exp(\hbar^{-1}p_1q_1+\CCh(\v)) \, d\mu(p) .
$$
\end{theorem}
\begin{proof}
Using the explicit heat kernel
$$
\exp\Bigl( \frac{t}{2} \frac{\p^2}{\p q^2} \Bigr) f(q)
= \int_{-\infty}^\infty f(p) \exp(-(p-q)^2/2t) \, \frac{dp}{\sqrt{2\pi t}} ,
$$
we see that, at a formal level,
$$
\exp(\Delta) f(\hbar,q_1,q_2,\dots) = \int_{\R^\infty}
\exp\Bigl( \sum_{n=1}^\infty \hbar^n \frac{\p}{\p p_{2n}} \Bigr) f(\hbar,p)
\prod_{n=1}^\infty \frac{e^{-(p_n-q_n)^2/2n\hbar^n}dp_n}{\sqrt{2\pi n\hbar^n}} .
$$
Using the formal integral $\hint_{\R^\infty}f\,d\mu(p)$, we may rewrite this
rigourously as
$$
\exp\Bigl( - \sum_{n=1}^\infty \frac{q_n^2}{2n\hbar^n} \Bigr)
\hint_{\R^\infty} \left\{ \exp\Bigl( \sum_{n=1}^\infty
\hbar^n \frac{\p}{\p p_{2n}} \Bigr) f(\hbar,p) \right\}
\exp\Bigl( \sum_{n=1}^\infty \frac{p_nq_n}{n\hbar^n} \Bigr) \, d\nu(p) .
$$
Integrating by parts, we obtain
\begin{multline*}
\exp\Bigl( - \sum_{n=1}^\infty \frac{q_n^2}{2n\hbar^n} \Bigr)
\hint_{\R^\infty} f(p)
\exp\Bigl( \sum_{n=1}^\infty \frac{p_nq_n}{n\hbar^n}
+ \sum_{\text{$n$ even}} \Bigl\{ \frac{p_n-q_n}{n\hbar^{n/2}} -
\frac{1}{2n} \Bigr\} \Bigr) \, d\nu(p) \\
\begin{split}
&= \exp\Bigl( - \sum_{n=1}^\infty \frac{q_n^2+q_{2n}}{2n\hbar^n} \Bigr)
\hint_{\R^\infty} f(p) \exp\Bigl( \sum_{n=1}^\infty \frac{p_nq_n}{n\hbar^n}
\Bigr)\, d\mu(p) \\
&= \Exp(-\hbar^{-1}h_2) \hint_{\R^\infty} f(p) \Exp(\hbar^{-1}p_1q_1) \,
d\mu .
\qed\end{split}
\end{multline*}
\def\qed{}
\end{proof}

Although it is possible that \ref{7.17}, the analogue for cyclic operads,
can be obtained from \ref{Wick} by the principle of stationary phase, we do
not know how to do this.

In the next section, we need the following consequence of \ref{Wick}.
\begin{corollary} \label{Wick:Det}
$$
- \hbar^{-1}e_2 + \CCh(\FF_{\Det}\a) = - \Log \hint_{\R^\infty}
\Exp \Bigl( \hbar^{-1}p_1q_1 - \CCh(\a) \Bigr) \, d\mu(p)
$$
\end{corollary}
\begin{proof}
By \ref{Det}, we have
$\Det^\Dual\cong\local_\Susp\o\local_\Sigma$. Applying \ref{tom}, we see
that
$$
\CCh(\FF_{\Det}\a) \cong
\CCh(\Sigma\Susp\MM\Susp^{-1}\Sigma^{-1}\a^*) \cong
- \tom \CCh(\MM\Susp\Sigma\a) .
$$
Since $\tom\*\Log=\Log\*\tom$, the result follows from \ref{Wick}.
\end{proof}

\newcommand{\kint}{\int^\#}

\section{Euler characteristics of moduli spaces of curves and
$\protect\CCh(\FF_{\Det}\Ass)$} \label{ribbon}

In this section, we apply the results of Section 8 to calculate
$\CCh(\FF_{\Det}\Ass)$ explicitly. Using the decompositions of moduli
spaces of curves found by Harer, Mumford, Penner and others, we obtain new
information on the Euler characteristics of these moduli spaces \ref{M/S}.

\begin{definition}
A ribbon graph is a graph $G$, each vertex of which has valence at least
$3$, together with a cyclic order on the set of flags $v$ making up each
vertex $v\in\VV(G)$. (This is what Penner \cite{Penner} calls a fat graph.)
\end{definition}

Equivalently, a ribbon structure on a graph $G$ is the same as an isotopy
class of embeddings of the CW complex $|G|$ into a compact oriented Riemann
surface $\Sigma(G)$ with boundary, such that
\begin{enumerate}
\item the intersection of the image of $|G|$ with the boundary
$\p\Sigma(G)$ is the set of endpoints of the legs of $|G|$;
\item the image of $|G|$ is a deformation retract of $\Sigma(G)$.
\end{enumerate}
The cyclic orders of the sets $v$ are then induced by the embedding
$|G|\hookrightarrow\Sigma(G)$ and the orientation of $\Sigma(G)$.

Denote by $\gamma(G)$ and $\nu(G)$ the genus and number of boundary
components of $\Sigma(G)$. Note that $2(\gamma(G)-1)+\nu(G)=g(G)-1$, where
$g(G)=\dim H_1(|G|)$.

Ribbon graphs are related to the operad $\Ass$ in the following way. Since
$\Ass\(n\)\cong\Ind_{\Z_n}^{\SS_n}(\k)$, it has a basis $\{e_\sigma\}$
labelled by the cyclic orders on $\{1,\dots,n\}$ (see \ref{ass}). It
follows that $\FF\Ass\(g,n\)$ has a natural basis $\{e_G\}$ labelled by all
ribbon $n$-graphs $G$ with $g(G)=g$.
\begin{proposition} \label{9.1.1}
If $\FF\Ass\(\gamma,\nu,n\)$ is the subcomplex spanned by ribbon graphs $G$
with $\gamma(G)=\gamma$ and $\nu(G)=\nu$, there is a splitting of chain
complexes
$$
\FF\Ass\(g,n\) = \bigoplus_{2(\gamma-1)+\nu=g-1} \FF\Ass\(\gamma,\nu,n\) .
$$
\end{proposition}
\begin{proof}
If $e$ is an edge of a ribbon graph $G$, there is a natural ribbon graph
structure on $G/e$. If $e$ is a loop, the vertex of $G/e$ corresponding to
$e$ has genus $1$, and since $\Ass\(g,n\)=0$ for $g>0$, we see that
$\Ass\(G/e\)=0$; thus, the term of the differential on $\FF\Ass$ associated
to this edge vanishes. On the other hand, if both ends of $e$ are distinct,
then $\gamma(G/e)=\gamma(G)$ and $\nu(G/e)=\nu(G)$, so that the
corresponding term of the differential preserves the splitting.
\end{proof}

We may also consider the cyclic operad $\Ass$ as a modular $\Det$-operad,
where $\Det$ is the cocycle $\Det(G)=\Det(H^1(G))$. (This uses the fact
that the cocycle $\Det$ is canonically trivial on trees; see \ref{DET}.) The
Feynman transform $\FF_{\Det}\Ass$ also has a basis labelled by ribbon
graphs, and we have a decomposition of $\FF_{\Det}\Ass$ into a sum of
subcomplexes $\FF_{\Det}\Ass\(\gamma,\nu,n\)$ similar to \ref{9.1.1}.

\subsection{$\FF\Ass$, $\FF_{\Det}\Ass$ and moduli spaces of curves} 
In this section, we relate the complexes $\FF\Ass\(\gamma,\nu,n\)$ and
$\FF_{\Det}\Ass\(\gamma,\nu,n\)$ to cell decompositions of moduli spaces of
punctured curves (which we learnt about from Penner). These moduli spaces
are differentiable orbifolds, by which we mean, by analogy to the algebraic
case \ref{DM}, a proper \'etale differentiable groupoid $\CG$.

A $\nu$-punctured curve is a pair $(\Sigma,A)$, where $\Sigma$ is a smooth
projective algebraic curve over $\C$, and $A\subset\Sigma$ is a finite
subset. An isomorphism of two punctured curves
$(\Sigma_1,A_1)\to(\Sigma_2,A_2)$ is an isomorphism $\Sigma_1\to\Sigma_2$
inducing a bijection $A_1\to A_2$. A frame $\lambda$ of a punctured curve
$(\Sigma,A)$ is an element of the circle bundle over $A$ whose fibre at
$z\in A$ is the quotient of $T_z\Sigma\backslash\{0\}$ by the dilatation
group $\R^\times_+$.

An $n$-framed punctured curve is an object
$(\Sigma,A,\lambda_1,\dots,\lambda_n)$, where $(\Sigma,A)$ is a punctured
curve, and $(\lambda_1,\dots,\lambda_n)$ are $n$ distinct frames in
$(\Sigma,A)$.  An isomorphism of $n$-framed punctured curves is defined in
the obvious way.

The groupoid of $n$-framed pointed curves of genus $\gamma$ such that
$|A|=\nu$, and their isomorphisms, represents a differentiable orbifold
$Q_{\gamma,\nu,n}$. In fact, using the method of level structures, this
groupoid is seen to be equivalent to a transformation groupoid (the
quotient, in the sense of orbifolds, of a space by a group action). Note
that for $n=0$,
$$
Q_{\gamma,\nu,0} = \CM_{\gamma,\nu}/\SS_\nu .
$$ 

A cellular decomposition of an orbifold $\CG$ is a cellular decomposition
of $\Ob(\CG)$ whose inverse images in $\Mor(\CG)$ under the \'etale maps
$s$ and $t:\Mor(\CG)\to\Ob(\CG)$ coincide. Associated to this decomposition
is a cochain complex $C^\bull(|\CG|)$, the invariants of the action of
$\Mor(\CG)$ on the cellular cochain complex of $C^\bull(\Ob(\CG))$. This
complex may be thought of as the cochain complex associated to the
decomposition of $|\CG|$ into orbicells: these are the image in $|\CG|$ of
cells in $\Ob(\CG)$, and are quotients of cells by finite groups.

In the cellular decompositions which we study, the cells will not
necessarily be relatively compact; thus, the cellular cochain complex
$C^\bull(\CG)$ calculates the cohomology with compact supports
$H^\bull_c(\CG)$; this is isomorphic to the cohomology with compact
supports of $|\CG|$, as long as we work over a field of characteristic
zero.

The following result was communicated to us by R. Penner. Only the cases
$g=0$ and $n=0$ may be found in the literature, and it is only the case
$n=0$ which we will need.
\begin{theorem} \label{penner}
For all $\gamma\ge1$, $\nu\ge0$ and $n\ge0$, there is an orbifold
$P_{\gamma,\nu,n}$, fibred over $Q_{\gamma,\nu,n}$ with fibres $\R_+^\nu$,
with a natural $\SS_n$-equivariant cell decomposition, and a natural
identification $C^\bull(|P_{\gamma,\nu,n}|)\cong\FF\Ass\(\gamma,\nu,n\)$.
\end{theorem}

Informally, $P_{\gamma,\nu,n}$ parametrizes Riemann surfaces of genus
$\gamma$ with $\nu$ boundary circles, together with $n$ numbered points on
the boundary; the circumferences of these circles are labelled by points in
the fibre $\R^\nu_+$.

\subsubsection{\boldmath$g=0$} \label{g=0}
In this case, only ribbon graphs with the topology of trees contribute to
$\FF\Ass\(0,n\)=\FF\Ass\(0,1,n\)\cong\Sigma\Susp\BB\Ass\(n\)$. Every ribbon
$n$-graph with the topology of a tree can be embedded into the plane,
inducing a cyclic order on the set $\{1,\dots,n\}$ of legs of the
graph. Thus, $\BB\Ass\(n\)$ splits into a sum of subcomplexes
$\BB\Ass\(n\)_\sigma$ labelled by cyclic orders $\sigma$.

On the other hand, $Q_{0,1,n}$ is the quotient of the configuration space
of $n$ distinct points in $S^1$ by rotations, and is the union of cells
$K_\sigma$ corresponding to cyclic orders $\sigma$ as above. Each cell
$K_\sigma$ may be identified with the interior of the Stasheff polytope
$K_{n+1}$ \cite{Stasheff}, and the cellular decomposition $K_\sigma$ of
\ref{penner} is (Poincar\'e) dual to the face decomposition of $K_{n+1}$,
since faces of $K_{n+1}$ correspond to planar trees.

\subsubsection{\boldmath$n=0$} \label{n=0}
For this case, we mention the references \cite{Penner} and
\cite{Kontsevich}. In the first of these, the orbifold $P_{\gamma,\nu}$
(decorated Teichm\"uller space) is constructed, and a cellular
decomposition given, whose cells are in bijection with isotopy classes of
so-called ``ideal cell decompositions'' of a fixed Riemann surface $\Sigma$
of genus $\gamma$ with $\nu$ punctures. As remarked on page 40 of Penner
\cite{Penner:euler}, the (Poincar\'e) dual of an ideal cell decomposition
is a ``spine'' on $\Sigma$, i.e.\ a graph $G$ in $\Sigma$ together with a
deformation retraction of $\Sigma$ to $G$. This shows that the cells in
Penner's decomposition of $P_{\gamma,\nu}$ are in bijection with ribbon
graphs $G$. As for the differential on $C^\bull(|P_{\gamma,\nu})$, we only
need the following result, which we establish explicitly.
\begin{proposition} \label{Penner}
$
H_c^i(\CM_{\gamma,\nu}/\SS_\nu) \cong
\Sigma^{1-2\gamma} H_{-i}(\FF_{\Det}\Ass\(\gamma,\nu,0\)) .
$
\end{proposition}
\begin{proof}
Let $\CL=\RR^*p_!\underline{\C}$ be the direct image with proper supports
along the fibres of the projection
$p:P_{\gamma,\nu,0}\to\CM_{\gamma,\nu}/\SS_\nu$. The graded sheaf $\CL$ is
an invertible local system on the orbifold $\CM_{\gamma,\nu}/\SS_\nu$,
concentrated in cohomological degree $\nu$. The Thom isomorphism shows that
$$
H^\bull_c(Q_{\gamma,\nu,0},\C) \cong
H^\bull_c(P_{\gamma,\nu,0},p^*\CL^{-1}) ;
$$
thus, the cohomology with compact supports of
$Q_{\gamma,\nu,0}\cong\CM_{\gamma,\nu}/\SS_\nu$ may be identified with the
homology of the cellular cochains on $P_{\gamma,\nu,0}$ with coefficients
in the graded local system $p^*\CL^{-1}$.

Let $G$ be a ribbon graph with no legs such that $2(\gamma-1)+\nu=g-1$, and
let $C_G$ be the corresponding cell of $P_{\gamma,\nu,0}$. We will prove
the  natural identification
\begin{equation} \label{twist-by-det}
H^\bull_c(C_G,p^*\CL^{-1}) \cong
\Sigma^{1-2\gamma} \bigl( \om(G) \o \Det(G)^{-1} \bigr) ,
\end{equation}
which identifies the complex of cellular cochains on $P_{\gamma,\nu,0}$
with coefficients in $p^*\CL^{-1}$ with
$\Sigma^{1-2\gamma}\FF_{\Det}\Ass\(\gamma,\nu,0\)$, proving the result.

To prove \ref{twist-by-det}, let $\Sigma(G)$ be the surface with boundary
constructed from the ribbon graph $G$. In Penner's theory, a point $z$ of
$C_G$ is represented by a hyperbolic metric on $\Sigma(G)$. The fibre
$p^{-1}(z)$ is homeomorphic to $\R^{\pi_0(\p\Sigma(G))}$, the coordinates
being logarithms of the lengths of the components of $\p\Sigma(G)$. Thus
$$
H^\bull_c(C_G,p^*\CL^{-1}) \cong H^\bull_c(C_G,\C) \o
H^\bull_c(p^{-1}(p(z)),\C)^{-1} .
$$
The graded vector space $H^\bull_c(C_G,\C)$ is naturally identified with
$\om(G)$, while the graded vector space $H^\bull_c(p^{-1}(z),\C)$ is
naturally identified with $\Det(\pi_0(\p\Sigma(G)))$. Thus,
\ref{twist-by-det} follows from the formula
$$
\Det(G)^{-1} \o \Det(\pi_0(\p\Sigma(G))) \cong \Sigma^{2\gamma-1} \C ,
$$
which we now prove.

Let $\bar\Sigma(G)$ be the compact surface obtained by gluing a disk $D_i$
along each component $S_i$ of $\p\Sigma(G)$. Define the determinant of a
finite-dimensional graded vector space $V_\bull$ to be
$$
\Det(V) = \Det(V_0) \o \Det(V_1)^{-1} \o \Det(V_2) \o \dots .
$$
If $M$ is a closed manifold, Poincar\'e duality gives a canonical
identification $\Det(H^\bull(M,\C))\cong\Sigma^{-e(M)}\C$, where $e(M)$ is
the Euler characteristic of $M$.

Consider the Mayer-Vietoris sequence for the decomposition of
$\bar\Sigma(G)$ as the union of $\Sigma(G)$ and $\coprod_iD_i$ along
$\coprod_iS_i$. From the multiplicativity of $\Det$ in long exact
sequences, we obtain the identification
$$
\Det(H^\bull(\bar\Sigma(G),\C)) \cong
\frac{\Det(H^\bull(\Sigma(G),\C))\o\Det(H^\bull(\coprod_iD_i,\C))}
{\Det(H^\bull(\coprod_iS_i,\C))} .
$$
By Poincar\'e duality, the denominator is trivial, while
$\Det(H^\bull(\bar\Sigma(G),\C))\cong\Sigma^{2(\gamma-1)}\C$. Since the
inclusion $|G|\hookrightarrow\Sigma(G)$ is a homotopy equivalence, we see
that
$$
\Det(H^\bull(\Sigma(G),\C)) \cong \Det(H^\bull(|G|,\C)) \cong
\Sigma^{-1}\Det(G)^{-1} ,
$$
completing the proof.
\end{proof}

\subsection{$\CCh(\FF_{\Det}\Ass)$: factorization of the integral}
By \ref{ass2} the characteristic of the cyclic operad $\Ass$ has a special
form: it is a sum of terms the $n$th of which depends on $p_n$ alone, in
the sense that
$$
h_1 + h_2 + \Ch(\Ass) =
- \sum_{n=1}^\infty \frac{\phi(n)}{n} \log(1-p_n) .
$$
The formula of \ref{Wick:Det} now gives
as
\begin{equation} \label{Wick:Ass}
- \hbar^{-1}e_2 + \CCh(\FF_{\Det}\Ass) = - \Log \hint_{\R^\infty} \Exp
\frac{1}{\hbar} \Bigl( p_1q_1 + \sum_{n=1}^\infty \frac{\phi(n)}{n}
\log(1-p_n) \Bigr) \, d\mu(p) .
\end{equation}
The special form of this formula will allow us to calculate it by a
separation of variables. We will then calculate these integrals separately
using the method of stationary phase.

\subsection{Stationary phase and Wick's theorem}
Let $f\in\Q\[x,\hbar\]$ be a power series of the form
$$
f = \frac{x^2}{2} + \sum_{2(g-1)+n>0} \frac{f_{g,n}\hbar^g x^n}{n!} .
$$
The exponential $\exp(-f/\hbar)$ has the form
\begin{equation} \label{>}
\sum_{k=-\infty}^\infty \sum_{\substack{\ell\ge0\\ 2k+\ell>0}} c_{k,\ell}
\, \hbar^k x^\ell \exp(-x^2/2\hbar) , \quad c_{k,\ell} \in \Q .
\end{equation}
\begin{definition}
Define the formal integral
$$
\log \kint \exp\bigl( \hbar^{-1}(x\xi-f(x,\hbar)) \bigr) \,
\frac{dx}{\sqrt{2\pi\hbar}} \in \hbar^{-1} \Q\[\xi,\hbar\]
$$
by the formula
$$
\log \sum_{k,\ell} c_{k,\ell} \, \hbar^k \int_{-\infty}^\infty x^k \,
e^{x\xi/\hbar-x^2/2\hbar} \, \frac{dx}{\sqrt{2\pi\hbar}} .
$$
This is well-defined by \ref{>}.
\end{definition}

\begin{remark}
The coefficients $F_{g,n}$ of the power series
$$
\log \kint \exp(\hbar^{-1}(x\xi-f(x,\hbar)) \, \frac{dx}{\sqrt{2\pi\hbar}}
= \frac{1}{\hbar} \Biggl( \frac{\xi^2}{2} + \sum_{2(g-1)+n>0}
\frac{F_{g,n}\hbar^g\xi^n}{n!} \Biggr)
$$
may be calculated by Wick's formula, mentioned in the introduction:
$$
F_{g,n} = \sum_{G\in\Gamma\(g,n\)} \frac{1}{|\Aut(G)|} \prod_{v\in\VV(G)}
f_{g(v),n(v)} .
$$
In particular, they are given by universal polynomials
$F_{g,n}\in\Q\[f_{g,n}\]$.
\end{remark}

We can now state the formula of stationary phase (a special case of Theorem
7.7.7 of H\"ormander \cite{Hormander}). The proof uses nothing more than
Taylor's theorem and integration by parts.
\begin{theorem} \label{stationary}
Let $\phi$ be a differentiable function of $x\in(a,b)$ and
$\hbar\in(-\eps,\eps)$, such that the function $x\mapsto\phi(x,0)$ has only
the single critical point $x=0$ in the interval $(a,b)$. Let $f$ be the
Taylor series of $\phi$ around $(x,\hbar)=(0,0)$, and suppose that
$f_{xx}(0,0)=1$. Let $u\in C^\infty_c(a,b)$ equal $1$ in a neighbourhood of
the critical point $0$. Then there is an asymptotic expansion
$$
\log \int_a^b u(x)\,\exp(\hbar^{-1}(x\xi-\phi(x,\hbar)) \,
\frac{dx}{\sqrt{2\pi\hbar}} \sim \log \kint
\exp(\hbar^{-1}(x\xi-f(x,\hbar)) \, \frac{dx}{\sqrt{2\pi\hbar}} .
$$
\end{theorem}

Let us give a simple application of Theorem \ref{stationary}.
\begin{proposition} \label{Stirling}
$$
\log \kint \exp \frac{1}{\hbar} \bigl( x\xi+x+\log(1-x) \bigr) \,
\frac{dx}{\sqrt{2\pi\hbar}} = \frac{1}{\hbar} \Bigl( \xi-\log(1+\xi) -
\hbar \, \log(1+\xi) + \sum_{g=2}^\infty \frac{\zeta(1-g)}{1-g} \hbar^g
\Bigr)
$$
\end{proposition}
\begin{proof}
The function $x+\log(1-x)$ has a unique critical point $0$ in the interval
$(-\infty,1)$. By the theorem of stationary phase, if $u\in
C^\infty_c(-\infty,1)$ equals $1$ in a neighbourhood of $0$,
$$
\log \int_{-\infty}^1 u(x) \, \exp \frac{1}{\hbar} \bigl( x\xi+x+\log(1-x)
\bigr) \, \frac{dx}{\sqrt{2\pi\hbar}} \sim \log \kint \exp \frac{1}{\hbar}
\bigl( x\xi+x+\log(1-x) \bigr) \, \frac{dx}{\sqrt{2\pi\hbar}} .
$$
In fact, we may take $u=1$, since the contribution of the integral away
from a neighbourhood of $0$ may be shown to vanish to infinite order by
repeated integration by parts.

Performing the changes of variables $u=(1-x)(1+\xi)/\hbar$ and
$s=\hbar^{-1}$, we see that
\begin{multline*}
\log \int_{-\infty}^1 \exp\bigl(\hbar^{-1}(x\xi+x+\log(1-x))\bigr) \,
\frac{dx}{\sqrt{2\pi\hbar}} \\
\begin{split}
& = \log s^{-s-\frac12} (1+\xi)^{-s-1} e^{s(1+\xi)} \int_0^\infty u^s
e^{-u} \, \frac{du}{\sqrt{2\pi}} \\ & = \log{} (2\pi)^{-1/2} s^{-s-\frac12}
(1+\xi)^{-s-1} e^{s(1+\xi)} \Gamma(s+1) \\ & = s(1+\xi) + \log{}
(2\pi)^{-1/2} s^{-s+\frac12} (1+\xi)^{-s-1} \Gamma(s) .
\end{split}
\end{multline*}
The proposition follows on inserting Stirling's formula
$$
\log\Gamma(s) \sim \sum_{k=1}^\infty \frac{\zeta(-k)}{-k} s^{-k}
+ \left(s-\frac{1}{2}\right) \log(s) - s + \frac{1}{2} \log(2\pi) ,
\quad s\to+\infty ,
$$
and replacing $s$ by $\hbar^{-1}$.
\end{proof}

\subsection{Application of stationary phase to $\CCh(\FF_{\Det}\Ass)$}
For $n\ge1$, introduce the Laurent polynomial
$$
\alpha_n(\hbar) = \frac{1}{n} \sum_{d|n} \frac{\phi(d)}{\hbar^{n/d}}
= \frac{1}{n\hbar^n} (1 + O(\hbar^{n/2}) ,
$$
and the formal integral
\begin{equation} \label{I_n}
I_n(q_n,\hbar) = \log \kint \exp \frac{1}{n\hbar^n} \bigl( p_nq_n + p_n +
n\hbar^n\alpha_n(\hbar) \log(1-p_n) \bigr) \,
\frac{dp_n}{e^{c(n)/2n}\sqrt{2\pi n\hbar^n}} .
\end{equation}
This power series may be transformed into one of the form which we
considered above, by making the changes of variables
$p_n=n^{1/2}\hbar^{\frac{n-1}2}x$ and $q_n=n^{1/2}\hbar^{\frac{n-1}2}\xi$,
which convert it to
\begin{multline*}
I_n(n^{1/2}\hbar^{\frac{n-1}2}\xi,\hbar) \\ = \log \kint \exp
\frac{1}{\hbar} \bigl( x\xi + n^{-1/2}\hbar^{-\frac{n+1}2}x +
\hbar\alpha_n(\hbar) \log(1-n^{1/2}\hbar^{\frac{n-1}2}x) \bigr) \,
\frac{dx}{e^{c(n)/2n}\sqrt{2\pi \hbar}} .
\end{multline*}
Strictly speaking, we need a slight generalization of the formal integral
$\log\kint\dots$, in which $f(x,\hbar)$ depends not on $\hbar$, but on
$\hbar^{1/2}$; this does not present any additional difficulties.

Using the power series $I_n$, we may rewrite \ref{Wick:Ass} as
\begin{multline} \label{drei}
- \hbar^{-1}e_2 + \CCh(\FF_{\Det}\Ass) \\
\begin{split}
& = - \sum_{n=1}^\infty \Log \kint \exp \frac{1}{n\hbar^n} \bigl( p_nq_n +
p_n + n\hbar^n\alpha_n(\hbar) \log(1-p_n) \bigr) \,
\frac{dp_n}{e^{c(n)/2n}\sqrt{2\pi n\hbar^n}} \\ & = - \sum_{n=1}^\infty
\sum_{\ell=1}^\infty \frac{\mu(\ell)}{\ell} I_n(q_{\ell n},\hbar^\ell).
\end{split}
\end{multline}
The following result resolves a problem left open in \cite{Kontsevich}.
\begin{theorem} \label{BIG}
Let $\Psi_n(\hbar)$, $n\ge1$, be the power series
$$
\Psi_n(\hbar) = \sum_{k=1}^\infty \frac{\zeta(-k)}{-k} \alpha_n^{-k}
+ (\alpha_n+1/2)\log(n\hbar^n\alpha_n) - \alpha_n + \frac{1}{n\hbar^n} -
\frac{c(n)}{2n} .
$$

\begin{enumerate}
\item The series
$$
\Psi(\hbar) = \sum_{n=1}^\infty \sum_{\ell=1}^\infty \frac{\mu(\ell)}{\ell}
\Psi_n(\hbar^\ell)
$$
is convergent; more precisely, $\Psi_n=O\bigl(\hbar^{\lceil n/6 \rceil}\bigr)$.
\item We have
$$
- \hbar^{-1}e_2 + \CCh(\FF_{\Det}\Ass) = \hbar^{-1}p_1 - \bigl(
\hbar^{-1} + 1 \bigr) \sum_{n=1}^\infty \frac{\phi(n)}{n} \log(1+p_n) -
\Psi(\hbar) .
$$
\end{enumerate}
\end{theorem}
\begin{proof}
Let $\beta_n=n\hbar^n\alpha_n-1$. Then
$\beta_n=c(n)\hbar^{n/2}+O(\hbar^{\lceil 2n/3 \rceil})$ and
$\beta_n^2=c(n)\hbar^{n/2}+O(\hbar^{\lceil 7n/6 \rceil})$, so that
\begin{align*}
(\alpha_n+1/2)\log(n\hbar^n\alpha_n) &= \frac{1}{n\hbar^n}
(1+\beta_n+\half n\hbar^n) \bigl( \beta_n - \beta_n^2/2 + \beta_n^3/3 + \dots
\bigr) \\
& = \frac{1}{n\hbar^n} \bigl( \beta_n + \beta_n^2/2 + O(\hbar^{3n/2})
\bigr) \\
& = \alpha_n - \frac{1}{n\hbar^n} + \frac{c(n)}{2n} + O\bigl(\hbar^{\lceil
n/6 \rceil} \bigr) .
\end{align*}
It follows that
$$
(\alpha_n+1/2)\log(n\hbar^n\alpha_n) - \alpha_n + \frac{1}{n\hbar^n} -
\frac{c(n)}{2n} = O\bigl(\hbar^{\lceil n/6 \rceil} \bigr) .
$$
On the other hand, the term proportional to $\alpha_n^{-k}$ in the
definition of $\Psi_n$ has the form $O(\hbar^{kn})$, and hence these terms
converge to a power series which does not contribute to the leading order
behaviour of $\Psi_n$. This completes the proof of part 1).

Generalizing \ref{Stirling}, which is the special case $n=1$, we see that
\begin{multline*}
I_n(q_n,\hbar) = \frac{q_n}{n\hbar^n} - (\alpha_n+1) \log(1+q_n) \\ +
\sum_{k=1}^\infty \frac{\zeta(-k)}{-k} \alpha_n^{-k} +
(\alpha_n+1/2)\log(n\hbar^n\alpha_n) - \alpha_n + \frac{1}{n\hbar^n} -
\frac{c(n)}{2n} ;
\end{multline*}
the only difference in the proof is that we make the substitutions
$u=(1-p_n)(1+q_n)/n\hbar^n$ and $s=\alpha_n(\hbar)$. Inserting this formula
into \ref{drei}, we see that
\begin{align*}
-\hbar^{-1}e_2+\CCh(\FF_{\Det}\Ass) =& - \sum_{n=1}^\infty
\sum_{\ell=1}^\infty \frac{\mu(\ell)q_{\ell n}}{\ell n\hbar^{\ell n}} -
\sum_{n=1}^\infty \sum_{\ell=1}^\infty \frac{\mu(\ell)}{\ell}
\alpha_n(\hbar^\ell) \log(1+q_{\ell n}) \\ & \quad - \sum_{n=1}^\infty
\sum_{\ell=1}^\infty \frac{\mu(\ell)}{\ell} \log(1+q_{\ell n}) -
\Psi(\hbar) .
\end{align*}
The definition of the M\"obius function shows that the first term equals
$$
- \sum_{k=1}^\infty \sum_{d|k} \frac{\mu(k/d)q_k}{k\hbar^k} =
- \frac{q_1}{\hbar} .
$$
Inserting the definition of $\alpha_n$ into the second term, we see that it
equals
$$
- \sum_{k=1}^\infty \log(1+q_k) \sum_{d|k}
\frac{\phi(d)}{k\hbar^{k/d}} \sum_{e|(k/d)} \mu(e)
= \sum_{k=1}^\infty \frac{\phi(k)}{k\hbar} \log(1+q_k) .
$$
The third term equals
$$
- \sum_{k=1}^\infty \frac{1}{k} \log(1+q_k) \sum_{d|k} d \mu(k/d)
= \sum_{k=1}^\infty \frac{\phi(k)}{k} \log(1+q_k) ,
$$
since by M\"obius inversion, $\sum_{d|k}\mu(k/d)d=\phi(k)$.
\end{proof}

By \ref{Penner}, this theorem implies the following purely topological
formula. (The change of sign comes from the fact that
$H^\bull_c(\CM_{g,n}/\SS_{g,n},\C)$ is an odd suspension of
$H_\bull(\FF_{\Det}\Ass)$.)
\begin{corollary} \label{M/S}
The power series $\Psi(\hbar)$ has the following topological
interpretation,
$$
\Psi(\hbar) = \sum_{g=2}^\infty \hbar^{g-1}
\sum_{2(\gamma-1)+\nu=g-1} e(|\CM_{\gamma,\nu}/\SS_\nu|) ,
$$
where $e(|\CM_{\gamma,\nu}/\SS_\nu|)$ is the Euler characteristic
of the topological space $|\CM_{\gamma,\nu}/\SS_\nu|$.
\end{corollary}

The first few terms of $\Psi(\hbar)$ are as follows:
$$
\Psi(\hbar) = 2\hbar + 2\hbar^2 + 4\hbar^3 + 2\hbar^4 + 6\hbar^5 + 6\hbar^6
+ 6\hbar^7 + \hbar^8 + O(\hbar^9) .
$$

\subsection{Remarks}
(a) Observe that the coefficients of $\hbar^n$, $n>0$, in \ref{BIG} are
constant (that is, independent of the power sums $p_i$). This is in
agreement with \ref{penner}, since the Euler characteristic of
$|Q_{\gamma,\nu,n}|$ vanishes if $n>0$, provided
$2(\gamma-1)+\nu>0$. Indeed, there is a circle action on
$Q_{\gamma,\nu,n}$, given by the formula
$$
(\Sigma,A,\lambda_1,\dots,\lambda_n) \mapsto
(\Sigma,A,e^{it}\lambda_1,\dots,e^{it}\lambda_n) .
$$
The isotropy groups of the induced circle action on $|Q_{\gamma,\nu,n}|$
are finite, since a punctured Riemann surface $(\Sigma,A)$ has finitely
many automorphisms fixing the punctures. Thus, all of the orbits of this
circle action are circles, and the Euler characteristic of
$|Q_{\gamma,\nu,n}|$ vanishes.

As expected, the coefficient of $\hbar^{-1}$ in $\CCh(\FF\Ass)$ is just
$-\tom\Ch(\Ass)$, consistent with the fact that
$\Cyc(\FF_{\Det}\Ass)=\Sigma\Susp\BB\Ass\simeq\Sigma\Susp\Ass$.

(b) The virtual Euler characteristic $\chi(\CG)$ of an orbifold $\CG$ may
defined using a cellular decomposition of $\CG$ as the sum
$$
\chi(\CG) = \sum_{\text{cells $\CU$ of $\Ob(\CG)$}}
\frac{(-1)^{\dim(\CU)}}{|\Aut(\CU)|} ,
$$
where $\Aut(\CU)$ is the group of morphisms of $\CG$ fixing $\CU$. In
general, the virtual characteristic is a rational number. It behaves well
under quotienting: if $G$ is a finite group acting on $\CG$,
$$
\chi(\CG/G) = \frac{\chi(\CG)}{|G|} .
$$

The virtual Euler characteristic of the orbifolds $\CM_{\gamma,1}$ was
shown by Harer and Zagier \cite{HZ} to equal $\zeta(1-2\gamma)$. Since the
virtual Euler characteristic is multiplicative for fibrations of orbifolds,
the fibrations
\begin{equation}\begin{diagram} \label{fibrations}
\node{\CM_{\gamma,\nu}} \arrow{e} \arrow{s} \node{\CM_{\gamma,\nu-1}} \\
\node{\CM_{\gamma,\nu}/\SS_\nu}
\end{diagram}\end{equation}
immediately imply formulas for the virtual Euler characteristics of
$\CM_{\gamma,\nu}$ and $\CM_{\gamma,\nu}/\SS_\nu$. However, on descent to
the coarse moduli spaces, the maps of \ref{fibrations} are not fibrations,
so there is no elementary relation between the Euler characteristics of the
topological spaces $|\CM_{\gamma,\nu}|$ and $|\CM_{\gamma,\nu}/\SS_\nu|$
for different $\nu$. (In fact, the Euler characteristics of the coarse
moduli spaces $|\CM_{\gamma,\nu}|$ and $|\CM_{\gamma,\nu}/\SS_\nu|$ are
unknown for $\nu>1$.)

Harer and Zagier also calculate the Euler characteristic of
$|\CM_{\gamma,1}|$ (page 482, \cite{HZ}), obtaining the formula
$$
\sum_{\gamma=1}^\infty e(|\CM_{\gamma,1}|) \hbar^{2\gamma-1}
= \sum_{n=1}^\infty \frac{\phi(n)}{n}
\sum_{\ell=1}^\infty \mu(\ell) \Psi_{n,\ell}(\hbar) ,
$$
where
$$
\Psi_{n,\ell}(\hbar) = \sum_{k=1}^\infty \zeta(-k) \alpha_{n,\ell}^{-k}
+ \alpha_{n,\ell}\log(n\hbar^n\alpha_n) + \frac{1}{n\hbar^n} -
\alpha_{n,\ell} .
$$
Here, $\alpha_{n,\ell}$ is the Laurent polynomial
$$
\alpha_{n,\ell}(\hbar) = \frac{1}{n} \sum_{d|n} \mu(d/(d,\ell))
\frac{\phi(n/d)}{\phi(\ell/(d,\ell))} \hbar^{-d} .
$$
There is a striking formal similarity between our formula \ref{M/S} for
$\Psi(\hbar)$ and this formula, representing the contribution of $\nu=1$.

(c) The idea of replacing an asymptotic integral (over an arbitrarily small
neighborhood of zero) by an equivalent integral over a fixed interval,
which is then explicitly evaluated through $\Gamma$-function and estimated
by Stirling formula, already arises in the calculation of the virtual Euler
characteristics of moduli spaces: see Harer and Zagier \cite{HZ}, Penner
\cite{Penner:euler}, Itzykson and Zuber \cite{IZ}, and Kontsevich
\cite{Kontsevich:KdV}.

\end{document}